\documentclass[12pt]{amsart}
\usepackage{amsmath,amsthm,amsfonts,amssymb}
\usepackage{graphicx}
\usepackage{mathrsfs}
\usepackage{enumitem}
\usepackage{etoolbox}
\usepackage{xcolor}
\apptocmd{\sloppy}{\hbadness 10000\relax}{}{}
\apptocmd{\sloppy}{\vbadness 10000\relax}{}{}
\usepackage[pdfpagelabels]{hyperref}
\usepackage[letterpaper,margin=1.1in]{geometry}

\numberwithin{equation}{section}
\theoremstyle{plain}

\newtheorem{theorem}{Theorem}[section]

\newtheorem{corollary}[theorem]{Corollary}

\newtheorem{lemma}[theorem]{Lemma}

\theoremstyle{definition}
\newtheorem{remark}[theorem]{Remark}
\newtheorem{definition}[theorem]{Definition}

\newcommand{\norm}[1]{\left\lVert#1\right\rVert}

\def\Xint#1{\mathchoice
{\XXint\displaystyle\textstyle{#1}}%
{\XXint\textstyle\scriptstyle{#1}}%
{\XXint\scriptstyle\scriptscriptstyle{#1}}%
{\XXint\scriptscriptstyle\scriptscriptstyle{#1}}%
\!\int}
\def\XXint#1#2#3{{\setbox0=\hbox{$#1{#2#3}{\int}$ }
\vcenter{\hbox{$#2#3$ }}\kern-.6\wd0}}

\def\dashint{\Xint-}

\newcommand{\diam}{\mathop\mathrm{diam}\nolimits}
\newcommand{\dist}{\mathop\mathrm{dist}\nolimits}

\newcommand{\loc}{\mathrm{loc}}

\DeclareMathOperator*{\esssup}{ess\,sup}

\begin{document}

\title{One-phase free-boundary problems with degeneracy}

\author{Sean McCurdy}

\address{Department of Mathematics\\ Carnegie Mellon University}
\email{seanmccu@andrew.cmu.edu}

\begin{abstract} In this paper, we study local minimizers of a degenerate version of the Alt-Caffarelli functional. Specifically, we consider local minimizers of the functional $J_{Q}(u, \Omega):= \int_{\Omega} |\nabla u|^2 + Q(x)^2\chi_{\{u>0\}}dx$ where $Q(x) = \text{dist}(x, \Gamma)^{\gamma}$ for $\gamma>0$ and $\Gamma$ a $C^{1, \alpha}$ submanifold of dimension $0 \le k \le n-1$.  We show that the free boundary may be decomposed into a rectifiable set, on which we prove upper Minkowski content estimates, and a degenerate cusp set about which little can be said in general with the current techniques.  Work in the theory of water waves and the Stokes wave serves as our inspiration, however the main thrust of this paper is to study the geometry of the free boundary for degenerate one-phase Bernoulli free-boundary problems in the context of local minimizers.

\end{abstract}

\maketitle

\tableofcontents

\renewcommand{\thepart}{\Roman{part}}

\section{Introduction}\label{sec:intro}

In the groundbreaking paper \cite{AltCaffarelli81} Alt and Caffarelli studied the existence and regularity of minimizers of the functional,
\begin{align}\label{alt-caffarelli functional}
J_Q(u, \Omega) := \int_{\Omega} |\nabla u |^2 + Q^2(x)\chi_{\{u > 0\}} dx
\end{align}
for $\Omega \subset \mathbb{R}^n$ an open, connected set with locally Lipschitz boundary $\partial \Omega$ such that the boundary satisfies  $\mathcal{H}^{n-1}(\partial \Omega) >0$. Minimization happens in the class, $$K_{u_0, \Omega} := \{u \in W^{1, 2}(\Omega) | u-u_0 \in W^{1,2}_0(\Omega) \},$$ for a $u_0 \in W^{1, 2}(\Omega)$ satisfying $u_0 \ge 0$, under the assumption that $Q$ was bounded and measurable.

In this paper, we study local minimizers of (\ref{alt-caffarelli functional}) where $Q(x) = \dist(x, \Gamma)^{\gamma}$ for any $\gamma > 0$, where $\Gamma$ is a $k$-dimensional $C^{1, \alpha}$ submanifold in $\mathbb{R}^n$, with $0 \le k \le n-1$.  In their original paper \cite{AltCaffarelli81}, the authors prove the regularity of the \textit{free boundary} $\partial \{u>0\} \cap \Omega$ 
under the assumption that the function $Q$ satisfies the following conditions:
\begin{enumerate}
\item $Q \in C^{0, \alpha}(\Omega)$
\item $0< Q_{\min} \le Q(x) \le Q_{\max} < \infty.$
\end{enumerate}
This assumption plays a crucial role in proving non-degeneracy of the minimizing function $u$ near the free boundary and the non-degeneracy of the free boundary itself.  For example, if $n = 2$ then near a free boundary point $x_0$, we can write $$u(x) = Q(x_0)\langle x-x_0, \vec \eta \rangle_+ + O(|x - x_0|)$$ for some unit vector $\vec \eta$. Thus, if $0< Q(x_0)$ the blow-up of $u$ at $x_0$ is a piece-wise linear function and hence the blow-up of $\partial \{u>0\}$ is flat. On the other hand, if $Q(x_0) = 0$ then we see that $u$ cannot have piece-wise a linear blow-up at $x_0$, and therefore all blow-ups of $\partial \{u>0\}$ are not flat.  We decompose $\partial \{u>0\}$ into the \textit{regular set}, $\text{reg}(\partial \{u>0\})$, where blow-ups of $\partial \{u>0\}$ are piece-wise linear and the \textit{singular set}, $\text{sing}(\partial \{u>0\})$, where blow-ups are not piece-wise linear. We shall refer to the case when $0< Q_{\min} \le Q$ as the \textit{non-degenerate case}.  

In higher dimensions, the non-degenerate case becomes more complicated, since there exist blow-ups which are not piece-wise linear.  However, the following results are known.  
\begin{theorem}\label{t:previous work}
Let $0< Q_{\min} \le Q \le Q_{\max}$ be $C^{0,\alpha}$ in $B_2(0) \subset \mathbb{R}^n.$  Suppose that $u$ is a local minimizer of $J_{Q}(\cdot, B_2(0))$. Then $\partial \{u>0\} \cap B_2(0)$ may be decomposed into 
\begin{align*}
    \partial \{u>0\} = \emph{reg}(\partial \{u>0\}) \cup \emph{sing}(\partial \{u>0 \}).
\end{align*}
\begin{itemize}
    \item[i.] (\cite{AltCaffarelli81}) The regular set, $\text{reg}(\partial \{u>0\})$, is relatively open and can be locally written as the graph of a $C^{1, \beta}$ function.  
    \item[ii.] (\cite{AltCaffarelli81}) In $n=2$, $\emph{sing}(\partial \{u > 0 \}) = \emptyset.$
    \item[iii.] (\cite{Weiss99})If $k^*$ is the first dimension that there exists a non-linear one-homogeneous minimizer of $J_{1}(\cdot, \mathbb{R}^n)$, then $\dim_{\mathcal{H}}(\emph{sing}(\partial \{u >0\})) \le n-k^*.$
    \item[iv.] (\cite{CaffarelliJerisonKenig04}) $k^* \ge 4$.
    \item[v.] (\cite{JerisonSavin15}) $k^* \ge 5 $.
    \item[vi.] (\cite{DeSilvaJerison09}) $k^* \le 7$
\end{itemize} 
Furthermore, if $Q$ is $C^{k, \alpha}$ (resp. smooth) for some $k \ge 1$ and $0< \alpha <1$, then $\emph{reg}(\partial \{u>0\})$ is locally the graph of a $C^{k+1, \beta}$ (resp. smooth) function \cite{AltCaffarelli81}.
\end{theorem}

The results of \cite{AltCaffarelli81} have inspired a countless papers and generalizations, including two-phase versions \cite{AltCaffarelliFriedman84}, fractional Laplacian versions \cite{CRS10} and \cite{DeSilvaRoquejoffre12}, $p$-Laplacian versions \cite{DanielliPetrosyan05}, and elliptic versions \cite{Valdinoci04}, and versions for almost-minimizer \cite{DavidToro15}\cite{DavidEngelsteinGarciaToro19},  just to name a few.  Much interest has centered around the exact value of $k^*$ in the non-degenerate case, and in \cite{EdelenEngelstein19}, Edelen and Engelstein prove that $\text{sing}(\partial \{u>0\})$ is $(n-k^*)$-rectifiable and satisfies certain upper Minkowski content estimates.  

However, virtually all work on the geometry of the related free boundaries has proceeded under the assumption of non-degeneracy on $Q$.  Indeed, to the author's knowledge \cite{AramaLeoni12} is the first instance in which the \textit{degenerate case}, i.e. $Q_{\min} = 0$, was considered. In \cite{AramaLeoni12}, Arama and Leoni investigate absolute minimizers and assume $n=2$, $Q(x_1, x_2) = \sqrt{(c - x_2)_+}$, $\Omega$ is a rectangle, and symmetric boundary conditions on $\partial \Omega$. Subsequent work for $n \ge 3$, $0 < \gamma$, and $Q(x', x_n) = (c - x_n)_+^{\gamma}$ has been carried out in \cite{GravinaLeoni18}\cite{GravinaLeoni19}, though also only for absolute minimizers under similar assumptions of symmetry. These investigations were carried out to investigate the theory of the water waves and the Stokes wave, which are usually merely critical points of $J_{Q}(\cdot, \Omega)$, in the context of minimizers. No work has been done for other functions $Q$.

Inspired by the work of \cite{AramaLeoni12} and \cite{GravinaLeoni18}, we study the fine-scale structure of $\text{sing}(\partial \{u>0 \})$ for local minimizers in the degenerate case for a natural class of $Q(x) = \dist(x, \Gamma)^\gamma$. In particular, we make no assumptions of symmetry and allow $\Gamma$ to be non-flat. Since our questions are essentially local and Alt and Caffarelli \cite{AltCaffarelli81}, Weiss \cite{Weiss99}, and Edelen and Engelstein \cite{EdelenEngelstein19} have proven detailed partial regularity results on $\text{reg}(\partial \{ u>0\})$ and $\text{sing}(\partial \{u>0\})$ when $Q_{\min} \ge c >0$, we restrict our investigation to questions on the infinitesimal structure of $\partial \{u >0\} \cap \{Q = 0\}.$

\subsection{Main results.}

Our most basic result allows us to decompose $\partial \{u>0\} \cap \Gamma$ into two sets which we must treat very differently.

\begin{lemma}\label{l:sigma and S decomp}
Let $0< \gamma < \infty$ and $0 \le k \le n-1$ be an integer.  Let $\Gamma$ be a $k$-dimsnesional $C^{1, \alpha}$ submanifold.  If $u$ is a local minimizer of (\ref{alt-caffarelli functional}) with $Q(x)=\dist(x, \Gamma)^{\gamma}$ then the following holds.
\begin{enumerate}
    \item For all $x \in \partial \{u>0\} \cap \Gamma$, the $Q$-density
    \begin{align}\label{e:Q-density def}
    \Phi(x, 0^+) := \lim_{r \rightarrow 0^+}\frac{1}{\omega_n r^{n+2\gamma}}\int_{B_r(x)}Q^{2}\chi_{\{u>0 \}}dx
\end{align}
exists.
    \item There exist constants $0< c(n, \gamma)< C(n, \gamma)< \int_{B_1(0)}|x_n|^{2\gamma}dx$ such that we may decompose $\partial \{u>0\} \cap \Gamma$ as
    \begin{align*}
    \partial \{u>0\} \cap \Gamma = \mathcal{S} \cup \Sigma,
\end{align*}
    where $\Sigma = (\partial \{u>0\} \cap \Gamma) \cap \{x : \Phi(x, 0^+) =0\}$ and $$\mathcal{S} = (\partial \{u>0\} \cap \Gamma) \cap \{x: \Phi(x, 0^+) \in [c(n, \gamma), C(n, \gamma)]\}.$$ Moreover, $\partial \{u>0\} \cap \Gamma \subset \emph{sing}(\partial \{u>0\})$, since all blow-ups at points in $\mathcal{S}$ are $(1+ \gamma)$-homogeneous and all blow-ups at points in $\Sigma$ vanish identically.
\end{enumerate}
\end{lemma}
The existence of the $Q$-density is proven in Lemma \ref{l:density}. The lower bound $0< c(n, \gamma)$ is proven in Lemma \ref{l:energy lower-bound}.  The upper bound is proved in Lemma \ref{Q-density upper bound}.

For lack of a better term, the set $\Sigma$ shall be refered to as is the set of \textit{degenerate} singularities, and the set $\mathcal{S}$ shall be call the set of \textit{non-degenerate} singularities.  

The main results of this paper concern the non-degenerate set $\mathcal{S}$.  Specifically, we prove upper Minkowski content estimates on the ``effective" strata of $\mathcal{S}^{j}_{\epsilon, r_0} \subset \mathcal{S}$ using the powerful techniques of \cite{NaberValtorta17}.  Roughly speaking, the strata $\mathcal{S}^{j}_{\epsilon, r_0}$ is the set of points $x \in \Gamma \cap \partial \{u>0\}$ such that for all $r_0 < r \le \dist(x, \partial \Omega)$ the function $u$ is ``$\epsilon$-far" in $B_r(x)$ from all homogeneous functions which are translation invariant along a $(k+1)$-dimensional linear subspace. See Section \ref{ss:quant strat} for details and rigorous definitions.

\begin{theorem}\label{t:main theorem}
Let $n \ge 2$, $k$ be an integer such that $0 \le k \le n-1$, and let $0< \gamma$.  Let $\Gamma \subset \mathbb{R}^n$ be a $k$-dimensional $(1, M)$-$C^{1, \alpha}$ submanifold such that $0 \in \Gamma$, and let $Q(x) = \dist(x, \Gamma)^{\gamma}$. 

If $u$ is an $\epsilon_0$-local minimizer of $J_Q(\cdot , B_2(0))$ in the class $K_{u_0, B_2(0)}$ and $\norm{\nabla u_0}^2_{L^2} \le 
\Lambda$ and $\sup_{\partial B_2(0)}u_0 = A<\infty.$  Then, for any $0< \epsilon, \rho$ and any radius $r$ such that $\rho \le r$,
\begin{align*}
\emph{Vol}(B_r(\mathcal{S}^j_{\epsilon, \rho} \cap B_{1}(0))) \le C(n, j, \epsilon, \gamma, \alpha, M, \epsilon_0, \Lambda, A)r^{n-j}.
\end{align*}
Furthermore, $\mathcal{S}^j_{\epsilon}= \bigcap_{0<\rho}\mathcal{S}^j_{\epsilon, \rho}$ is countably $j$-rectifiable.  Consequently, for all $0 \le j \le k$ the strata $\mathcal{S}^{j}:= \bigcup_{0<\epsilon}\mathcal{S}^{j}_{\epsilon}$ are countably $j$-rectifiable.
\end{theorem}

Thanks to Lemma \ref{l:sigma and S decomp}(2), we are able to prove the following containment result.

\begin{lemma}\label{l:e reg}($\epsilon$-containment)
Under the hypotheses of Theorem \ref{t:main theorem}, $\mathcal{S} \cap \overline{B_1(0)}$ is closed and satisfies the following containment relationship. If $k = n-1$, then there exists an $0< \epsilon(n, \gamma, \alpha, M, \Lambda, A, 
\epsilon_0)$ such that $\mathcal{S} \subset \mathcal{S}^{n-2}_{\epsilon}.$
\end{lemma}

Theorem \ref{t:main theorem} and Lemma \ref{l:e reg} immediately imply the following corollary.
\begin{corollary}\label{c: finite upperminkowski bounds}
Under the hypotheses of Theorem \ref{t:main theorem}, if $k =n-1$, $\overline{dim_{\mathcal{M}}}(\mathcal{S})\le n-2$ and 
\begin{align*}
    \mathcal{H}^{n-2}(\mathcal{S} \cap B_1(0)) \le \mathcal{M}^{*, n-2}(\mathcal{S} \cap B_1(0)) \le C(n, \gamma, \alpha, M, \Lambda, A, \epsilon_0).
\end{align*}  
\end{corollary}
We note that when $k < n-1$ similar estimates are obtainable, but are meaningless, as there exists simple examples in which $\mathcal{S} = \Gamma \cap B_1(0)$. When $k = n-1$, these result say that the non-degenerate singular set $\mathcal{S}$ cannot by too spread out, since it must ``sit in space" like an $(n-2)$-dimensional submanifold, but also cannot be too concentrated at any scale, either.

\subsubsection{Cusps}

A central concern in the degenerate case is the formation of degenerate singularities, i.e. cusps.  Since the classical estimates in \cite{AltCaffarelli81} which are essential to regularity become vacuous in regions where $Q$ vanishes, the standard analytic techniques of establishing weak geometric regularity (i.e., interior ball conditions) do not work near $\Gamma$. New-- or at least different-- ideas are needed.

The potential development of cusps leads to some notable differences between the degenerate case and the non-degenerate case.  See, Lemma  \ref{l:positivity set convergence} and Remark \ref{r:positivity set convergence}. 

In this paper, we apply analytic techniques to prove the following preliminary result on the degenerate singular set $\Sigma$.
\begin{lemma}\label{l:finite perimeter 0}
If $u$ is a local minimizer of $J_Q(\cdot, B_2(0))$, then $\partial \{u>0 \} \cap B_1(0)$ is a set of finite perimeter.  In particular, the set $\Sigma = \Gamma \cap \partial \{u>0\} \setminus \mathcal{S}$ has $\mathcal{H}^{n-1}$-measure zero.
\end{lemma}

\begin{remark}
Subsequent to writing this paper, the author and Lisa Naples have proven some general conditions under which $\Sigma = \emptyset$. In particular, when $n \ge 2$, $0 \le k \le n-1$ is an integer, $0< \gamma$, and $\Gamma \subset \mathbb{R}^n$ a \emph{flat} $k$-dimensional submanifold, then $\Sigma = \emptyset$ for all $\epsilon_0$-local minimizers of $J_Q(\cdot , B_2(0))$ for $Q(x) = \dist(x, \Gamma)^{\gamma}$. See \cite{MccurdyNaples22}.

Notwithstanding the \textit{post facto} non-existence of cusps under these circumstances, the techniques in this paper may prove useful in other circumstances when cusps cannot be eliminated.
\end{remark}

\subsection{Strategy and Organization}

The overall strategy of the paper is to employ the tools and techniques of \cite{NaberValtorta17} and to prove the density ``gap" in Lemma \ref{l:sigma and S decomp}(2).  This density ``gap" allows us to prove the $\epsilon$-containment results in Corollary \ref{l:e reg}. The key ingredients leading to these results are: non-degeneracy of $u$ and local Lipschitz estimates depending upon $Q(x) = \dist(x, \Gamma)^{\gamma}$. With these results, one is able to show that at non-degenerate singularities the the Weiss $(1+\gamma)$-density is almost monotone and points have $(1+\gamma)$-homogeneous blow-ups. The degree of homogeneity roughly follows from the fact that $Q$ is $\gamma$-homogeneous and $Q \sim |\nabla u|$.  

In Section \ref{s:preliminaries}, we introduce the quantitative stratification from \cite{CheegerNaber13} which is necessary to the machinery of \cite{NaberValtorta17}, as well as state some basic estimates on $C^{1, \alpha}$ submanifolds that are necessary to establish the almost-monotonicity formula for the Weiss $(1+\gamma)$-density on $\partial \{u>0\} \cap \Gamma$. Section \ref{s:existence} is dedicated to recapping the results of Alt and Caffarelli \cite{AltCaffarelli81} and Arama and Leoni \cite{AramaLeoni12}.  In particular, these comparison techniques prove the non-degeneracy and local Lipschitz bounds adjusted to the function $Q$. 

Section \ref{s:Weiss density} proves the almost-monotonicity of the Weiss $(\gamma +1)$-density for local minimizers and variational solutions. With the almost-monotonicity formula, we then prove strong compactness and the existence of $(1+ \gamma)$-homogeneous blow-ups in Section \ref{s:companctness}.  Note that the degeneracy of $Q$ plays an important role in the compactness results (see Remark \ref{r:positivity set convergence}).  In Section \ref{s:density}, we prove the lower bounds in the density ``gap" in the definition of $\mathcal{S}$, see Lemma \ref{l:energy lower-bound}.  

Section \ref{s:proof of main theorem} is devoted to carefully following the argument of \cite{EdelenEngelstein19}, which applied the techniques of \cite{NaberValtorta17} to the non-degenerate case to prove Theorem \ref{t:main theorem}, noting all the changes necessary to adapt them to the degenerate case. In Section \ref{s:e reg}, we prove the containment results which prove Lemma \ref{l:e reg}. Finally, in Section \ref{s:finite perimeter}, we prove that the positivity set $\{u >0\}$ is a set of finite perimeter, and hence $\Sigma$ has $\mathcal{H}^{n-1}$-measure zero.

\subsection{Acknowledgements}
The author acknowledges the Center for Nonlinear Analysis at Carnegie Mellon University for its support.  In particular, the author thanks Giovanni Leoni and Irene Fonseca for their invaluable generosity, patience, and guidance.  In addition, the author thanks the reviewer for many helpful and excellent comments.

\section{Preliminaries}\label{s:preliminaries}

In this section, we discuss some basic definitions and elementary results.

\subsection{$C^{1,\alpha}$ geometry}\label{s:manifold geometry}
 In this subsection, we record some elementary observations on the submanifolds $\Gamma$ which are necessary to prove the almost-monotonicity results in Lemma \ref{t:Weiss almost monotonicity}.
\begin{definition}
Let $n \in \mathbb{N}$ and $0\le k \le n-1$.  A set $\Gamma$ is locally a $k$-dimensional $C^{1, \alpha}$ submanifold if for every $x \in \Gamma$, there is a radius $0< r_x$ such that
\begin{align*}
\Gamma \cap B_{r_x}(x) = \text{graph}_{T_x\Gamma}(f_x)
\end{align*}
for a some function, $f_x \in C^{1, \alpha}(\mathbb{R}^k; \mathbb{R}^{n-k})$. We will use $[f]_{\alpha, B_{r_x}(x)}$ to denote the H\"older seminorm of $Df$ in $B_{r_x}(x) \subset \mathbb{R}^k$.  That is,
\begin{align*}
[f]_{\alpha, B_{r_x}(x)}:= \sup_{z, y \in B^k_{r_x}(x), z \not = y} \frac{|Df(y) - Df(z)|}{|y-z|^{\alpha}}.
\end{align*}

We shall call $\Gamma$ a $(r, M)$-$C^{1, \alpha}$ submanifold of dimension $0 \le k \le n-1$, if for every $x \in \Gamma$, we may take $r_x = r$ and $ [f]_{\alpha, B_{r}(x)} \le M.$  For $\Gamma \subset \mathbb{R}^n$ a $(r, M)$-$C^{1, \alpha}$ submanifold of dimension $0\le k \le n-1$, we denote
\begin{align*}
[\Gamma]_{\alpha} := \sup \{[f_x]_{\alpha, B_{r}(x)} : x \in \Gamma \} \le M.
\end{align*}
\end{definition}

\begin{lemma}\label{l:Gamma containment}
Let $\Gamma \subset \mathbb{R}^{n}$ be a $(1, M)$-$C^{1, \alpha}$ submanifold of dimension $0\le k \le n-1$.  Then, for any $x \in \Gamma$ with defining function $f: B_1(x) \cap (x + T_x\Gamma) \rightarrow \mathbb{R}^{n-k}$, for all $y \in B_1(x) \cap (x + T_x\Gamma)$, 
\begin{align*}
|(y, f(y))-(y, 0)| & = [\Gamma]_{\alpha}|y|^{1+ \alpha}. 
\end{align*}
\end{lemma}

\begin{proof}
By translation and rotation, we may assume that $x = 0$ and $\frac{\partial}{\partial x_j}f(0) = 0$ for $1 \le j \le k$.  That is, $T_x\Gamma = \mathbb{R}^k \hookrightarrow \mathbb{R}^n$.  We calculate,  
\begin{align*}
|(y, f(y))-(y, 0)| & = |\int_0^y Df(z)\cdot \frac{y}{|y|}dz|\\
& \le |y| \max_{z \in \overline{0y}}||Df(z)|| = [\Gamma]_{\alpha}|y|^{1+ \alpha}. 
\end{align*}
\end{proof}

\begin{remark}\label{r:dist gradient}
The function, $\dist(\cdot, \Gamma): \mathbb{R}^n \rightarrow \mathbb{R}$, is a Lipschitz function with Lipschitz constant $1.$  By Rademacher's theorem, $\nabla(\dist(\cdot, \Gamma))$ exists $\mathcal{H}^n$-a.e..  Furthermore, for $x \in \mathbb{R}^n$ such that $\nabla(\dist(x, \Gamma))$ exists, there exists a \textit{unique} minimizing point $y \in \Gamma$ such that $\dist(x, \Gamma) = |x-y|$ and
\begin{align*}
\nabla(\dist(x, \Gamma)) = \frac{x- y}{|x- y|}.
\end{align*}
\end{remark}

\begin{definition}
We define the function, $\pi_{\Gamma}: \Omega \rightarrow \Gamma$, as follows. For $x$ such that $\dist(\cdot, \Gamma)$ is differentiable at $x$, we define $\pi_{\Gamma}(x):= y$, where $y\in \Gamma$ is the unique minimizing point such that $\dist(x, \Gamma) = |x-y|$.  By Remark \ref{r:dist gradient}, this is sufficient to define the function $\mathcal{H}^n$-a.e..
\end{definition}

\begin{remark}
Let $\Gamma \subset \mathbb{R}^{n}$ be a locally $C^{1, \alpha}$ submanifold of dimension $0\le k \le n-1$.  Then, for every $x \in \Omega$ such that $\pi_{\Gamma}(x)$ is defined, 
\begin{align*}
x-\pi_{\Gamma}(x) \perp T_{\pi_{\Gamma}(x)}\Gamma.
\end{align*}
\end{remark}

\begin{lemma}\label{l:normal drift}
Let $\Gamma \subset \mathbb{R}^{n}$ be an $(r_0, M)$-$C^{1, \alpha}$ submanifold of dimension $0\le k \le n-1$.  Let $x_0 \in \Gamma$, then for every $x \in B_{r_0}(x_0)$ such that $\pi_{\Gamma}(x)$ is defined, 
\begin{align*}
(\pi_{\Gamma}(x) - x_0) \cdot \frac{x - \pi_{\Gamma}(x)}{|x - \pi_{\Gamma}(x)|} \le 8[\Gamma]_{\alpha}r_0^{1 + \alpha}.
\end{align*}
\end{lemma}

\begin{proof}
We begin by choosing coordinates so that $\Gamma$ is a graph of $f:\mathbb{R}^k \rightarrow \mathbb{R}^{n-k}$ over $x_0 + T_{x_0}(\Gamma)$ where $x_0 = 0$ and $\frac{\partial}{\partial x_i}f(0) = 0$ for all $i = 1, 2, ..., k.$  We shall use the notation, $\vec \eta = \frac{x - \pi_{\Gamma}(x)}{|x - \pi_{\Gamma}(x)|}$. 
Let $y \in \mathbb{R}^k$ be such that $\pi_{\Gamma}(x) = (y, f(y)) \in \mathbb{R}^k \times \mathbb{R}^{n-k}.$  We decompose 
\begin{align*}
\pi_{\Gamma}(x) = \vec a + \vec b
\end{align*}
where $\vec a = (y, 0)$ and $\vec b = (0, f(y))$. Thus $\vec b \cdot \vec \eta  |\vec b| \le [\Gamma]_{\alpha}|y|^{1 + \alpha}$. Next, we choose a vector $\vec \eta_0 \in N_{0}\Gamma$, the normal bundle to $\Gamma$ at $0$, such that $|\vec \eta - \vec \eta_{0}| \le [\Gamma]|y^\alpha.$
\begin{align*}
\vec a \cdot \vec \eta & \le \vec a \cdot (\vec \eta + \vec \eta_0 - \vec \eta_0)\\
& \le \vec a \cdot \vec \eta_0 + \vec a \cdot \vec (\vec \eta - \vec \eta_0)\\
& \le [\Gamma]_{\alpha}|y|^{1+\alpha}.
\end{align*}
Since $|y| \le 2r_0$ and $0<\alpha \le 1,$ we have the desired inequality.
\end{proof}

\begin{lemma}(Manifold compactness)\label{l:manifold convergence}
Let $\Gamma_i$ be a sequence of $k$-dimensional $(1, M)$-$C^{1, \alpha}$ submanifolds satisfying $0 \in \Gamma_i$ for all $i \in \mathbb{N}.$  There is a $k$-dimensional $(1, M)$-$C^{1, \alpha}$ submanifold $\Gamma \subset \mathbb{R}^n$ with $0 \in \Gamma$ such that for every $0< R< \diam(B)$, 
\begin{align*}
\Gamma_j \cap B_R(0) \rightarrow \Gamma \cap B_R(0)
\end{align*}
in the Hausdorff metric on compact subsets.
\end{lemma}

This lemma is is a consequence of the compactness of the Grassmanian and the compactness of $C^{1, \alpha}$ functions with bounded seminorm $[f]_{\alpha, B_1(0)}$. The rest is left to the reader.

\begin{lemma}\label{l:distance estimate}
Let $\Gamma$ be a $k$-dimensional $(1, \frac{1}{4})$-$C^{1, \alpha}$ submanifold satisfying $0 \in \Gamma$.  Then, there is a constant $C_1< \infty$, such that 
\begin{align*}
    \mathcal{H}^{n}(B_{2^{-i}}(\Gamma) \cap B_{1}(0)) \le C_1 2^{-i(n-k)}.
\end{align*}
\end{lemma}

\begin{proof}
This fact holds more generally for Ahlfors regular sets.  We argue below for $\Gamma$, specifically.  We begin with an initial estimate, $$\mathcal{H}^{n}(B_{1}(\Gamma) \setminus B_{\frac{1}{2}}(\Gamma) \cap B_1(0)) \le \omega_{k}\omega_{n-k}.$$ Note that since $\Gamma$ is assumed to be $(1, \frac{1}{4})$-$C^{1, \alpha}$ submanifold, by Lemma \ref{l:Gamma containment}, we have that,
\begin{align*}
    (B_1(0) \cap B_{1}(\Gamma) \setminus B_{\frac{1}{2}}(\Gamma)) \subset B_{1}(0) \setminus B_{\frac{1}{4}}(T_0\Gamma),
\end{align*}
and hence $\mathcal{H}^n(B_1(0) \cap B_{1}(\Gamma) \setminus B_{\frac{1}{2}}(\Gamma)) \le \omega_{k}\omega_{n-k}.$

We now iterate this estimate at dyadic scales.  Note that since $\mathcal{H}^{k}(B_2(0) \cap \Gamma) \le 2^k\omega_{k}\sqrt{1 + [\Gamma]_{\alpha}} = C \omega_{k}$, an $r$-net in $\Gamma$ consists of at most $C(n, k) \omega_{k}r^{-k}$ points.  Furthermore, if $x \in B_1(0) \cap \{x: \dist(x, \Gamma)^{\gamma} \in (2^{-i},2^{-i+1}] \},$ then there must exist a $y \in \Gamma \cap B_{1+ 2^{-i+2}}(0)$ such that $x \in B_{2^{-i+1}}(y)$.  Additionally, we note that if $x \in \Gamma$, then $[\Gamma^{x, r}]_{\alpha} \le [\Gamma]_{\alpha}r^{\alpha}$.  Therefore, if we take an $2^{-i}$-net in $\Gamma \cap B_2(0)$ and apply our previous result to $B_{1}(0)$ and $\Gamma^{x_i, 2^{-i+1}}$, we obtain
\begin{align*}
    \mathcal{H}^{n}(B_1(0) \cap B_{2^{-i+1}}(\Gamma) \setminus B_{2^{-i}}(\Gamma)) & \le C(n, k)\omega_{k}2^{ik}(\omega_{k}\omega_{n-k})2^{(-i+1)n}\\
    & \le C(n, k)\omega_k^2\omega_{n-k}2^{n}2^{-i(n-k)}.
\end{align*}
This proves the estimate with $C_1 = C(n, k)\omega_k^2\omega_{n-k}2^{n}$.
\end{proof}

\subsection{Minimizers and local minimizers}

\begin{definition}\label{AC def} Let $\Omega$ be an open set with Lipschitz boundary satisfying $\mathcal{H}^{n-1}(\partial \Omega) >0$. Suppose that $Q \in C^{0,\alpha}(\Omega)$ such that $0 \le Q$.  A function $u$ is a \textit{minimizer} of 
\begin{align}
J_Q(u, \Omega) := \int_{\Omega} |\nabla u |^2 + Q^2(x)\chi_{\{u > 0\}} dx,
\end{align}
in the class $K_{u_0, \Omega}: = \{u \in W^{1, 2}(\Omega) | u-u_0 \in W^{1,2}_0(\Omega) \}$ for a $u_0 \in W^{1, 2}(\Omega)$ satisfying $u_0 \ge 0$, if for every other function $v \in K_{u_0, \Omega},$
\begin{align*}
J_{Q}(u, \Omega) \le J_{Q}(v, \Omega).
\end{align*} 

For $0< \epsilon_0$, a function $u$ is called an $\epsilon_0$-\textit{local minimizer} of $J_Q(\cdot, \Omega)$ if 
\begin{align*}
J_{Q}(u, \Omega) \le J_{Q}(v, \Omega)
\end{align*} 
for every $v \in K_{u, \Omega}$ satisfying 
\begin{align}
\norm{\nabla (u-v)}^2_{L^2(\Omega)} + \norm{\chi_{\{u>0\}} - \chi_{\{v>0\}}}_{L^1(\Omega)} < \epsilon_0. 
\end{align}
When we do not need to quantify such things, $\epsilon_0$-local minimizers will simply be called local minimizers. Clearly, minimizers are local minimizers for all $0<\epsilon_0$. And, for all $\Omega' \subset \Omega$, local minimizers in $K_{u, \Omega}$ are local minimizers in $K_{u, \Omega'}.$ We shall often speak of $u$ as a local minimizer without reference to the class $K_{u, \Omega}$.
\end{definition}

\begin{theorem}(\cite{AltCaffarelli81} Theorem 1.3)
Let $\Omega$ be an open set with Lipschitz boundary satisfying $\mathcal{H}^{n-1}(\partial \Omega) >0$. Suppose that $Q \in C^{0,\alpha}(\Omega)$ such that $0 \le Q$, and let $u_0 \in W^{1, 2}(\Omega)$ be non-negative satisfying $J_{Q}(u_0, \Omega)< \infty.$  Then, minimizers of the functional $J_{Q}(\cdot, \Omega)$ in the class $K_{u_0, \Omega}$ exist.
\end{theorem}

\subsection{Quantitative Stratification}\label{ss:quant strat}
Stratification is used in dimension-reduction arguments of the kind introduced by Federer and Almgren. It was applied to the non-degenerate case by Weiss in \cite{Weiss99} to show Theorem \ref{t:previous work}(\textit{iii.}), among other results.  These dimension-reduction techniques have been augmented into powerful ``effective" versions, first by Cheeger and Naber \cite{CheegerNaber13} to study manifolds with Ricci curvature bounded below, and then greatly strengthened by Naber and Valtorta \cite{NaberValtorta17} in the context of minimal surfaces. Following the deep analogy between minimal surfaces and local minimizers of (\ref{alt-caffarelli functional}) established by \cite{AltCaffarelli81}, this ``effective" version was used by Edelen and Engelstein \cite{EdelenEngelstein19} to address local minimizers of the non-degenerate case, as well as two-phase and vector versions.

The key to the improvement introduced by Cheeger and Naber is the quantitative control introduced in the form of a quantitative stratification.  See Definition \ref{def quant stratification}, below.

\begin{definition}(Symmetric functions and rescalings)
Given an integer $0 \le j \le n-2$, a function $u \in C^{0}(\mathbb{R}^n)$ is called $j$-\textit{symmetric} if $u$ is homogeneous and there is a linear subspace $V$ with $\text{dim}(V) = j$ such that for all $y \in V$,
\begin{align*}
    u(x) = u(x + y).
\end{align*}
For a function $u \in W^{1, 2}_{loc}(B_2(0)) \cap C^{0}(B_2(0))$ and a non-trivial ball $\overline{B_r(x)} \subset B_2(0)$, we define the rescaling of $u$ at $x$ at scale $0<r$ by,
\begin{align}\label{def T rescalings}
    T_{x, r}u(y) : = \frac{u(x + ry)}{\left(\int_{\partial B_1(0)}u(x + ry)^2d\sigma(y) \right)^{\frac{1}{2}}}.
\end{align}
 We shall say that $u$ is $(j, \epsilon)$-\textit{symmetric} in $B_r(x)$ if there exists a non-trivial $j$-symmetric function $\phi$ such that 
 \begin{align}
 \norm{T_{x, r}u - T_{0, 1}\phi}_{L^{2}(B_1(0))} \le \epsilon.
 \end{align}
\end{definition}

\begin{definition}\label{def quant stratification}(Quantitative stratification)
Let $n \ge 2$ and $0 \le k \le n-1$ be integers.  Let $0< \gamma< \infty$, $\Gamma$ be a $k$-dimensional $C^{1, \alpha}$ submanifold, and $Q(x)= \dist(x, \Gamma)^{\gamma}$. If $u$ is a local minimizer of \eqref{alt-caffarelli functional}, we make the following definitions.

For each $\epsilon > 0, r_0> 0$, and integer $0 \le j \le n-1$, we define the $(j, \epsilon, r_0)$-strata $\mathcal{S}_{\epsilon, r_0}^{j}$ as follows,
\begin{align*}
    \mathcal{S}_{\epsilon, r_0}^{j} := \{x \in \Gamma \cap \partial \{u>0\} : \text{$u$ is \emph{not} } (k+1, \epsilon)\text{-symmetric in $B_{r}(x)$ for all } r_0 \le r \le \dist(x, \partial \Omega)\}.
\end{align*}
That is, $x \in \mathcal{S}_{\epsilon, r_0}^{j}$ if and only if $\norm{T_{x, r}u - T_{0, 1}\phi}_{L^2(B_1(0))} \ge \epsilon$ for all non-trivial $(j+1)$-symmetric functions $\phi$. We shall use the notation $\mathcal{S}_{\epsilon}^j := \mathcal{S}^j_{\epsilon, 0}$. Observe that $\mathcal{S}^j := \bigcup_{0<\epsilon} \mathcal{S}^{j}_{\epsilon}$ is the traditional, qualitative $j$-stratum. Where we compare the strata of different functions, we shall use the notation $\mathcal{S}_{\epsilon, r_0}^{j}(u)$.
\end{definition}

\begin{remark}
The quantitative strata defined above are closed under $L^2_{loc}$ convergence of the underlying functions. In addition, they enjoy the following properties.
\begin{enumerate}
    \item $\mathcal{S}^{0} \subset \mathcal{S}^{1} \subset ... \mathcal{S}^{n-2} \subset \mathcal{S}^{n-1}$ and $\mathcal{S} \subset \mathcal{S}^{n-2}.$
    \item For all $\delta< \epsilon$ and $r<R$, $\mathcal{S}^j_{\epsilon, r} \subset \mathcal{S}^j_{\delta, R}$.  Furthermore, for integers $0 \le j < k \le n-1$, $\mathcal{S}^{j}_{\epsilon, r} \subset \mathcal{S}^k_{\epsilon, r}.$
\end{enumerate}
We note that by definition, $\partial \{u>0\} \cap \text{reg}(\partial \{u>0\}) \subset \mathcal{S}^{n-1}.$
\end{remark}

\section{Minimizers and Integral Average Growth Estimates}\label{s:existence}

In this section, we recall several important results for local minimizers of $J_{Q}(\cdot, \Omega)$. The main results in this section are the integral average growth estimates first established in \cite{AltCaffarelli81} which are crucial to establishing local Lipschitz estimates and non-degeneracy of the functions.  The results of this section are limited to local minimizers.

\subsection{Basic properties}

\begin{lemma}\label{outer variation}(Outer Variation, \cite{AltCaffarelli81} Lemma 2.2 and Lemma 2.3)
Suppose that $Q \in C^{0,\alpha}(\Omega)$ such that $0 \le Q$.  Minimizers and local minimizers of $J_{Q}(\cdot, \Omega)$ in the class $K_{u_0, \Omega}$ are subharmonic ($\Delta u \ge 0$ in $\Omega$) in the sense of distributions.  That is, for all $\phi \ge 0$, $\phi \in C^{\infty}_c(\Omega)$,
\begin{align*}
\int_{\Omega}\nabla u \cdot \nabla \phi dx \le 0.
\end{align*}
Moreover, $u$ is harmonic in the sense of distributions in the interior of $\{u>0\}$. And, we may choose a representative of $u \in W^{1, 2}(\Omega)$ which is defined point-wise by subharmonicity
$$
u(x) = \lim_{r\rightarrow 0}\dashint_{B_r(x_0)}udx.
$$
In this case, $u \in L^{\infty}_{loc}(\Omega)$ satisfying
\begin{align*}
    0 \le u(x) \le \esssup_{\Omega}u_0. 
\end{align*}
\end{lemma}

We now record the Noether Equations (inner variation) associated to being a local minimizer of $J_{Q}(\cdot, \Omega)$. See \cite{Velichkov19} Lemma 9.5 for the details of the calculation when $\gamma = 0$.

\begin{lemma}\label{inner variation}(Inner variation)
Let $0\le k\le n-1$ be integers. Let $0< M<\infty$, $0<\gamma$, $\Gamma$ be a $k$-dimensional $(1, M)$-$C^{1, \alpha}$ submanifold, and $Q(x) = \dist(x, \Gamma)^{\gamma}$.  For $u$ a local minimizer of $J_{Q}(\cdot, \Omega)$, and any $\phi \in C^{1}_0(\Omega; \mathbb{R}^n),$
\begin{align}\nonumber \label{e:Noether equations}
 0& = \int_{\Omega}\left(|\nabla u|^2 + Q^{2} \chi_{\{u >0\}}\right)\emph{div}(\phi)dx - \int_{\Omega}2\nabla u D\phi \nabla u  + \chi_{\{u>0\}} \nabla (Q^2(x)) \cdot \phi dx \\
& = \int_{\Omega}\left(|\nabla u|^2 + Q^{2} \chi_{\{u >0\}}\right)\emph{div}(\phi)dx\\
\nonumber & \qquad  \qquad  - \int_{\Omega} 2\nabla u D\phi \nabla u + 2 \gamma \dist(x, \Gamma)^{2\gamma-1} \chi_{\{u>0\}} \nabla (\dist(x, \Gamma))\cdot \phi dx.
\end{align}
\end{lemma}


\subsection{Integral Average Growth Estimates}

The techniques used in \cite{AltCaffarelli81} to establish the non-degeneracy of a local minimizer $u$ rely upon comparing $u$ with two other functions.  First, we will need to compare $u$ with the harmonic extension of $u$ in a ball $B_r(0).$ Second, we will need to compare $u$ with a function $w = \min\{u, v\}$ in $B_r(0)$ for

\begin{align*}
v(x) = \left(\sup_{y \in B_{r \sqrt{s}}(0)} \{u(y)\}\right) \max \biggl\{1 - \frac{|x|^{2-n} - r^{2-n}}{(sr)^{2-n} - r^{2-n}} , 0 \biggr\}.
\end{align*} 

\begin{remark}
Since $\norm{\nabla v}^2_{L^2} \le C(s, n) \sup_{y \in B_{r \sqrt{s}}(0)} \{u^2(y)\} r^{n-2}$ and harmonic functions are energy minimizers, for every $\epsilon_0$-local minimizer $u$, there is a uniform scale 
\begin{align} \label{e:standard scale}
    r_0= r_0(n, \sup_{\partial \Omega} u_0, \norm{\nabla u}_{L^2}, \epsilon_0)
\end{align}
at which we can apply these arguments. We shall refer to this scale $r_0$ as the \textit{standard scale}.
\end{remark}

\begin{theorem}(\cite{AltCaffarelli81} Lemma 3.2)
Let $n \ge 2$, and let $u$ be a $\epsilon_0$-local minimizer of $J_{Q}(\cdot, \Omega)$.  There is a constant, $C_{max}(n)>0$ such that for every $0< r< r_0$ and every ball $B_r(x) \subset \Omega$ if
\begin{align*}
\dashint_{\partial B_r(x)}ud\sigma  > C_{max}r \cdot \max_{y \in B_r(x)} Q(y)
\end{align*}
then $u > 0$ in $B_r(x)$.  In particular, if $\mathcal{H}^n(\{u=0\} \cap B_r(x))>0$ then 
\begin{align*}
\dashint_{\partial B_r(x)}ud\sigma < C_{max}r \cdot \max_{y \in B_r(x)} Q(y).
\end{align*}
\end{theorem}

We now argue that functions which satisfy the preceding properties are Lipschitz continuous.

\begin{lemma}\label{l:gradient upper bound} (cf. \cite{AltCaffarelli81} Corollary 3.3)
Let $n \ge 2$, $\Omega \subset \mathbb{R}^n$ be an open bounded set, and $Q \in C^{0, \alpha}(\Omega)$ be a non-negative function. Let $u$ be an $\epsilon_0$-local minimizer $J_{Q}(\cdot, \Omega)$.  Let $r_0$ be the standard scale.  Let $\Omega_{r} = \{x \in \Omega : \dist(x, \partial \Omega) > r\}.$  Then, $u$ enjoys the following properties.
\begin{enumerate}
    \item $u$ is harmonic in the classical sense in any $B_r(x) \subset \Omega$ for which $\mathcal{H}^n(\{u=0 \} \cap B_r(x)) = 0$.
    \item $\{u> 0 \}$ is open.
    \item $u$ is locally Lipschitz, satisfying
\begin{align*}
    \norm{\nabla u}_{L^{\infty}(\Omega_{2r_0})} \le \max\{ C(n)\frac{\norm{u}_{L^1(\Omega_{r_0})}}{r_0^{n+1}}, C(n)C_{max} \max_{y \in B_{r_0}(\partial \{u>0\})} \{Q(y) \}\}.
\end{align*}
\end{enumerate}
\end{lemma}

\begin{proof}
First, observe that since $u$ is a local minimizer, it must be that for all $0< r \le r_0$ $J_{Q}(u, B_r(x_0)) \le J_{Q}(h_u^{x_0, r}, B_{r}(x_0))$, where $h_u^{x_0, r}$ is the harmonic extension of $u$ in $B_r(x_0).$. Thus, using the orthogonality of harmonic functions,
\begin{align*}
    \int_{B_r(x_0)}|\nabla u - \nabla h^{x_0, r}_u|^2 dx \le \int_{B_r(x_0)}Q^2(x)\chi_{\{u=0\}}(x)dx.
\end{align*}
Therefore, if $\mathcal{H}^n(\{u=0 \} \cap B_r(x)) = 0$, then $u = h^{x_0, r}_u$. This proves (1).

To see that $\{u>0\}$ is open, let $x_0 \in \Omega \cap \partial \{u > 0\}$.  By (1) and the Maximum Principle, there must be an $0<\rho$ such that $\mathcal{H}^n(\{u=0 \} \cap B_r(x)) > 0$ for all $0< r<\rho$. Therefore, we may let $r \rightarrow 0$, we may invoke Lemma \ref{outer variation} and the pointwise definition of $u(x_0)$ to obtain
\begin{align*}
u(x_0) & \le \dashint_{\partial B_r(x_0)}ud\mathcal{H}^{n-1} \le C_{max} r \max_{y \in B_r(x_0)}\{ Q(y)\} \le C_{max}r \max_{y \in \Omega}\{Q(y)\} \rightarrow 0,
\end{align*}
where the last line follows from the fact that $|Q|$ is bounded since $Q \in C^{0, \alpha}(\Omega)$ and $\Omega$ is bounded.  Thus, $\{u > 0\} \cap \partial \{u>0\} = \emptyset$ and $\{u> 0 \}$ is open.

The local Lipschitz bound follows from standard harmonic estimates.  We break into two cases.  For any $x_0 \in \{u>0\} \cap \Omega_{2r_0} \setminus B_{r_0/2}(\partial \{u>0\}),$ we estimate
\begin{align*}
    |\nabla u(x_0)| & \le C_n \frac{2^{n+1}}{r_0^{n+1}}\int_{B_{r_0/2}(x_0)} udx\\
    & \le C_n \frac{1}{r_0^{n+1}}\norm{u}_{L^1(\Omega_{r_0})}.
\end{align*}

Now, consider $x_0 \in \{u>0 \} \cap \Omega_{2r_0} \cap B_{r_0/2}(\partial \{u>0\}).$ Let $0< \delta \le r_0/2$ and $x \in \partial \{u > 0\}$ be such that $|x - x_0| = \dist(x_0, \partial \{u>0\}) = \delta.$  Then, $u$ is harmonic in $B_{\delta}(x_0)$, and  
\begin{align*}
|\nabla u(x_0)| & \le C_n \frac{1}{\delta^{n+1}}\int_{B_{\delta}(x_0)} udx\\
& \le C_n \frac{1}{\delta^{n+1}}\int_{B_{2\delta}(x)} udx\\
& \le C(n) \frac{1}{\delta^{n+1}} \int_{0}^{2\delta} \int_{\partial B_{r}(x)} ud\mathcal{H}^{n-1}dr\\
& \le C(n) \frac{1}{\delta^{n+1}} \int_{0}^{2\delta} C_{max}\omega_{n-1}r^{n-1} r \max_{y \in B_r(x)}\{ Q(y)\}  dr\\
& \le C(n) \frac{C_{max}}{\delta^{n+1}}  \max_{y \in B_{2\delta}(x)}\{ Q(y)\} \int_{0}^{2\delta} \omega_{n-1}r^{n}dr\\
& \le C(n)C_{max} \max_{y \in B_{r_0}(\partial \{u>0 \})}\{ Q(y)\}.
\end{align*}
\end{proof}

\begin{corollary}\label{c:local Lipschitz bound}
Let $0 \le k \le n-2$, and let $\Gamma$ be a $k$-dimensional $(1, M)$-$C^{1, \alpha}$ submanifold such that $\Gamma \cap \Omega \not = \emptyset$.  Let $0< \gamma$, $Q(x) = \dist(x, \Gamma)^{\gamma}$, and $u$ a $\epsilon_0$-local minimizer of $J_{Q}(\cdot, \Omega)$ with standard scale $0< r_0$.  If $x_0 \in \partial \{u>0\}$, then for all $0 < r< \frac{1}{2}r_0$ 
\begin{align*}
||\nabla u||_{L^{\infty}(B_r(x_0))} \le C(n) r \max_{w \in B_{2r}(x_0)}\{ Q(w)\}.
\end{align*}
And, for all $x \in \{u>0\} \cap B_r(x_0)$
\begin{align*}
 |\nabla u(x)| \le C(n)2^{\gamma}\max\{\dist(x, \partial \{u>0\}), \dist(x, \Gamma)\}^{\gamma}.
\end{align*}
\end{corollary}

\begin{proof}
This corollary follows immediately from our choice of $Q$ and noting that in the penultimate line of Lemma \ref{l:gradient upper bound}, we have for $\delta = \dist(x_0, \partial \{u>0\})$
\begin{align*}
\max_{y \in B_{2\delta}(x)}\{ Q(y)\} \le 2^{\gamma}\max\{\dist(x, \partial \{u>0\}), \dist(x, \Gamma)\}^{\gamma}. 
\end{align*}
\end{proof}

\begin{lemma}\label{l:qmin ineq}(\cite{AltCaffarelli81} Lemma 3.4)
Let $u$ be an $\epsilon_0$-local minimizer of $J_Q(\cdot, \Omega)$.  Let $s \in (0, 1)$ be fixed.  Then, for all $0< r \le \frac{1}{2}r_0$ and all $B_{2r}(x_0) \subset \Omega$, there is a constant $C_{\min} = C(n, s)$ such that if,
\begin{align*}
\dashint_{\partial B_r(x_0)}u d\sigma \le r C_{\min} \min_{y \in B_{sr}(x_0)}\{ Q(y)\},
\end{align*}
then $u = 0$ on $B_{sr}(x_0).$ In particular, if $\{u>0 \} \cap B_{sr}(x_0) \not = \emptyset$, then 
\begin{align*}
\dashint_{\partial B_r(x_0)}u d\sigma > r C_{\min} \min_{y \in B_{sr}(x_0)}\{ Q(y)\}.
\end{align*}
\end{lemma}

\begin{corollary}\label{c:u growth}(Non-degeneracy of functions)
Let $\Gamma$ be a $(1, \frac{1}{4})$-$C^{1, \alpha}$ submanifold of dimension $0\le k \le n-1$ such that $0 \in \Gamma$.  Let $Q(x) = \dist(x, \Gamma)^{\gamma}.$ Let $u$ be a $\epsilon_0$-local minimizer of $J_{Q}(\cdot, B_2(0)).$ Then for every $0< r< r_0$ and $x \in \{u>0\} \cap \left( B_{r_0}(\partial \{u>0\})\setminus B_{r}(\partial \{u>0\}) \right).$  Then there is a constant $0<C(n, \gamma)< 1$ such that 
\begin{align*}
    u(x) \ge C(n, \gamma) r \max\{Q(x), (r/2)^{\gamma}\}.
\end{align*}
\end{corollary}

\begin{proof}
We break the proof into two cases.  First, suppose that $B_{\frac{1}{2}r}(x) \cap \Gamma = \emptyset.$ Then since $u> 0$ on $B_{\frac{1}{4}r}(x)$, then by Lemma \ref{l:qmin ineq} with $s = 1/2$, we have that
\begin{align*}
u(x) = \dashint_{\partial B_{r/2}(x)}u d\sigma \ge C_{\min}\frac{1}{2}r \min_{w \in B_{\frac{1}{4}r}(x)}\{Q(w)\} \ge C_{\min}(\frac{1}{2})^{2\gamma+1}rQ(x).
\end{align*}
If, on the other hand, $B_{\frac{1}{2}r}(0) \cap \Gamma \not = \emptyset$, then let $y \in \Gamma \cap B_{\frac{1}{2}r}(x)$.  By Lemma \ref{l:normal drift}, we have that $\Gamma \cap B_{\frac{1}{2}r}(x) \subset B_{\frac{1}{4}r^{1+\alpha}}(L)$ for some affine $k$-plane $L$ intersecting $B_{\frac{1}{2}r}(x)$.  Thus, we may find a point $y' \in \partial B_{\frac{1}{2}r}(x)$ such that $B_{\frac{1}{4}r}(y') \cap \Gamma = \emptyset.$  Note that we may choose $y'$ such that $Q(y') = \max_{w \in \partial B_{\frac{1}{2}r}(x)}\{Q(w)\}$. Thus,
\begin{align*}
u(x) & = \dashint_{\partial B_{\frac{1}{2}r}(x)}u d\sigma\\
& \ge (\omega_{n-1}(r/2)^{n-1})^{-1}\int_{\partial B_{\frac{1}{2}r}(x) \cap B_{\frac{1}{16}r}(y')}u d\sigma\\
& \ge (\omega_{n-1}(r/2)^{n-1})^{-1}\int_{\partial B_{\frac{1}{2}r}(x) \cap B_{\frac{1}{16}r}(y')}\left(\dashint_{\partial B_{r/16}(z)} u(w)d\sigma(w) \right) d\sigma(z)\\
& \ge C_{\min}\frac{r}{16^{n-1}}\min_{w \in B_{\frac{1}{8}r}(y')}\{ Q(w)\} \ge C_{\min}16^{1-n}r \left(\frac{r}{8}\right)^{\gamma}.
\end{align*}
\end{proof}

\begin{lemma}\label{l:interior balls}(Interior Balls)
Let $\Gamma$ be a $(1, \frac{1}{4})$-$C^{1, \alpha}$ submanifold of dimension $0\le k \le n-1$ such that $0 \in \Gamma$.  Let $Q(x) = \dist(x, \Gamma)^{\gamma}.$ Suppose that $u$ is an $\epsilon_0$-local minimizer of $J_{Q}(\cdot, B_2(0))$ and $x \in \partial \{u>0 \}$. Then for every $0< r\le r_0$, if $B_{r}(x) \subset B_1(0) \setminus \Gamma$, then there exists a point $y \in \{u>0 \} \cap \partial B_{\frac{1}{2}r}(x)$ and a constant $$0< c(n, \min_{B_{\frac{1}{2}r}(x)}\{Q(w)\}, \max_{B_{\frac{3}{4}r}(x)}\{Q(w)\}) <\frac{1}{2}$$ such that for all $z \in B_{\frac{c}{2}r}(y)$,
\begin{align*}
    u(z) \ge C_{min}(n) r \min_{B_{\frac{1}{2}r}(x)}\{Q(w)\},
\end{align*}
where $C_{\min}(n)$ is the constant from Lemma \ref{l:qmin ineq}.
\end{lemma}

\begin{proof}
By Lemma \ref{l:qmin ineq} with $s = \frac{1}{2}$,
\begin{align*}
    \dashint_{\partial B_{\frac{1}{2}r}(x)}ud\sigma \ge C_{\min}(n)\frac{1}{2}r \min_{w \in B_{\frac{1}{2}r}(x)}\{ Q(w)\}.
\end{align*}
Thus, there must exist a point $y \in \partial B_{\frac{1}{2}r}(x) \cap \{u>0\}$ such that $$u(y) \ge C_{\min}\frac{1}{2}r\min_{w \in B_{\frac{1}{2}r}(x)}\{ Q(w)\}.$$  By Lemma \ref{c:local Lipschitz bound}, $\text{Lip}(u|_{B_{\frac{3}{4}r}}(x)) \le C(n) r \max_{B_{\frac{3}{4}r}(x)}\{Q(w)\}$, it must be the case that for $$c = \frac{C_{min}\min_{w \in B_{\frac{1}{2}r}(x)}\{ Q(w)\}}{4C(n)\max_{B_{\frac{3}{4}r}(x)}\{Q(w)\}}$$ on $B_{cr}(z)$ we have that $u(z) \ge C_{min}(n)\frac{1}{4}r \min_{B_{\frac{1}{2}r}(x)}\{Q(w)\}$, as claimed.
\end{proof}

\section{Weiss Density}\label{s:Weiss density}

The Weiss densities were first introduced in \cite{Weiss99} Theorem 3.1 for critical points of more general, two-phase versions of non-degenerate $J_{Q}(\cdot, \Omega)$ functionals.  Since then Weiss-type densities, also called boundary adjusted energy functionals, have become an important tool in Bernoulli-type free boundary problems.  Their key property is that they give a monotonicity formula for critical points of $J_Q(\cdot, \Omega)$, thereby extending the analogy between minimal surfaces and Bernoulli free-boundary problems established by Alt and Caffarelli \cite{AltCaffarelli81} for regularity to the analysis of singularity.

In \cite{WeissVarvaruca11} Theorem 3.5, the authors extend this monotonicity formula to so-called ``weak solutions" to a free boundary problem associated to $J_Q(\cdot, \Omega)$ for the case $Q(x', x_n) = \sqrt{|x_n|}$. \cite{AramaLeoni12} applied this monotonicity formula to local minimizers of $J_{Q}(\cdot, \Omega)$ for $n=2$, $\gamma = 1/2$, and $\Gamma$ flat.  The main result of this section, Lemma \ref{t:Weiss almost monotonicity} below, extends this monotonicity formula to the cases $0<\gamma<\infty$ and $Q(x) = \dist(x, \Gamma)^\gamma$ where $\Gamma$ is a $k$-dimensional $(1, M)$-$C^{1, \alpha}$ submanifold. The calculation largely follows \cite{WeissVarvaruca11} Theorem 3.5, with the additional error carried around by fact that we now allow $\Gamma$ to be non-flat in a controlled way.

At the end of this Section, we give a brief discussion of Lemma \ref{t:Weiss almost monotonicity} and Corollary \ref{c:almost monotonicity} in comparison with the standard results for $\gamma =0$. See Remark \ref{r: gamma 0 non 0}.

For the remainder of this section, we assume $n \ge 2$ and $\Gamma$ a $(1, M)$-$C^{1, \alpha}$ submanifold of dimension $k \in \{0, 1, ..., n-1 \}$ such that $\Gamma \cap \Omega \not = \emptyset$.  Furthermore, assume that $0< \gamma < \infty$ and $Q(x) = \dist(x, \Gamma)^\gamma$.   

\begin{definition}(Weiss Density)
For any open set $\Omega \subset \mathbb{R}^n$ and any function $u \in W^{1, 2}(\Omega) \cap C(\Omega)$, we define the Weiss $(1 + \gamma)$-density with respect to $J_{\Omega}(\cdot, \Omega)$ at a point $x_0 \in \Omega$ and scale $0<r$ such that $\overline{B_r(x_0)} \subset \Omega$ as follows. 
\begin{align*}
W_{\gamma + 1}(x_0, r, u, \Gamma) = \frac{1}{r^{n-2+2(\gamma +1)}} \int_{B_r(x_0)}|\nabla u|^2 + Q^2(x)\chi_{\{u>0\}}dx - \frac{\gamma + 1}{r^{n-1+2(\gamma + 1)}}\int_{\partial B_r(x_0)} u^2 d\sigma.
\end{align*}
\end{definition}
\begin{remark}\label{r:Weiss invariances}
The Weiss density is invariant in the following senses.  For $u_{x, r}(y)= \frac{1}{r^{1+\gamma}}u(ry + x),$ and $\Gamma^{x, r} = \frac{1}{r}(\Gamma - x_0),$
\begin{align*}
W_{\gamma +1}(0, 1, u_{x, r}, \Gamma^{x, r}) = W_{\gamma + 1}(x, r, u, \Gamma)
\end{align*}
\end{remark}

\begin{lemma}\label{t:Weiss almost monotonicity}
Let $u$ is local minimizer of $J_{Q}(\cdot, \Omega)$ and $x_0 \in \Gamma \cap \partial \{u>0\}$. For almost every $0< r \le \frac{1}{2}\dist(x, \Omega^c)$ we have the following formula.
\begin{align}\label{Weiss derivative}\nonumber
& \frac{d}{dr}W_{\gamma +1}(x_0, r, u, \Gamma) \\
& \qquad \qquad =  \frac{2}{r^{n + 2\gamma}}\int_{\partial B_r(x_0)}(\nabla u \cdot \eta - \frac{\gamma + 1}{r}u)^2  + \\\nonumber
& \qquad \qquad \qquad \frac{1}{r^{n + 1 + 2\gamma}}2\gamma \int_{B_r(x_0)}\dist(x, \Gamma)^{2\gamma-1} \chi_{\{u>0\}} \frac{x - \pi_{\Gamma}(x)}{|x - \pi_{\Gamma}(x)|}\cdot (\pi_{\Gamma}(x) - x_0) dx.
\end{align}
Moreover,
\begin{align}\label{e:W derivative inequality}
\frac{d}{dr}W_{\gamma +1}(x_0, r, u, \Gamma) & \ge \frac{2}{r^{n+2\gamma}}\int_{\partial B_r(x_0)}\left(\nabla u \cdot \eta - \frac{\gamma +1}{r} u \right)^2d\sigma - 16 \gamma [\Gamma]_{\alpha}r^{\alpha - 1}.
\end{align}
\end{lemma}

\begin{proof}
Without loss of generality, by translation, we assume that $x_0 = 0$.  We begin by defining,
\begin{align*}
U(r) = \frac{1}{r^{n +2\gamma}}\int_{B_r(0)}|\nabla u|^2 + Q^2\chi_{\{u>0\}}dx, \qquad V(r) = \frac{1}{r^{n +1+2\gamma}}\int_{\partial B_r(0)}u^2 d\sigma
\end{align*}
so that $W_{\gamma + 1}(0, r, u, \Gamma) = U(r) - (\gamma + 1)V(r).$ We note that,
\begin{align}\label{e: U'}
U'(r) = & -\frac{n + 2\gamma}{r^{n + 1 +2\gamma}}\int_{B_r(0)}|\nabla u|^2 + Q^2\chi_{\{u>0\}}dx +  \frac{1}{r^{n+2\gamma}}\int_{\partial B_r(0)}|\nabla u|^2 + Q^2\chi_{\{u>0\}}dx.\\ \label{e: V'}
V'(r) = & -\frac{2 + 2\gamma}{r^{n + 2 +2\gamma}}\int_{\partial B_r(0)}u^2d \sigma + \frac{1}{r^{n +1 +2\gamma}}\int_{\partial B_r(0)}2 u \nabla u \cdot \eta d\sigma,
\end{align}
where $\eta$ is the unit outer normal.

Our first task is to use the inner variation formula to rewrite $U'(r)$. Let $\phi_{\tau} \in C^{\infty}_{c}(B_r(0))$ be a function which satisfies the following conditions:  
\begin{align*}
    \phi_{\tau}= 1 \text{ in } B_{r-\tau}(0), \qquad \nabla \phi_{\tau}(x) = -\frac{1}{\tau}\frac{x}{|x|} + o(\tau) \text{ in } B_r(0) \setminus B_{r-\tau}(0).
\end{align*}
For example, a suitable modification of $\phi(x) = \max\{0, \min \{1, \frac{r-|x|}{\tau}\} \}$ suffices. Define a vector field $\xi_{\tau}(x) = x \phi_{\tau}(x)$.  Observe that
\begin{align*}
    \text{div}(\xi_\tau(x)) & = n \phi_{\tau}(x) + x \cdot \nabla \phi_{\tau}(x)\\
    D \xi_{\tau}(x) & = \phi_\tau(x) Id_n + x \otimes \nabla \phi_\tau(x).
\end{align*}
Therefore, plugging $\xi_{\tau}$ into the Noether equations of Lemma \ref{inner variation} have 
\begin{align*}
0 = & \int_{\Omega}|\nabla u|^2 + Q^{2} \chi_{\{u >0\}} (n \phi_{\tau}(x) + x \nabla \phi_{\tau}(x))dx\\
& \qquad - 2\int_{\Omega} |\nabla u|^2 \phi_{\tau}(x) + \left(\nabla u \cdot x\right) (\nabla u \cdot \nabla \phi_{\tau}(x))dx\\
& \qquad + 2 \gamma \int_{\Omega} \dist(x, \Gamma)^{2\gamma-1} \chi_{\{u>0\}} \nabla (\dist(x, \Gamma))\cdot \xi_{\tau}(x) dx.
\end{align*}
Letting $\tau \rightarrow 0$, we have
\begin{align*}
0 = &n\int_{B_r(0)}|\nabla u|^2 + Q^{2} \chi_{\{u >0\}}dx  - r\int_{\partial B_r(0)}|\nabla u|^2 + Q^{2} \chi_{\{u >0\}}d\sigma\\
& \qquad - 2\int_{B_r(0)} |\nabla u|^2dx  + 2r\int_{\partial B_r(0)} (\nabla u \cdot \eta)^2 d\sigma\\
& \qquad  + 2 \gamma \int_{B_r(0)} \dist(x, \Gamma)^{2\gamma-1} \chi_{\{u>0\}} \frac{x - \pi_{\Gamma}(x)}{|x - \pi_{\Gamma}(x)|}\cdot x dx.
\end{align*}
Splitting the last term, by writing $x = x - \pi_{\Gamma}(x) + \pi_{\Gamma}(x)$, we obtain
\begin{align*}
& 2 \gamma \int_{B_r(0)} \dist(x, \Gamma)^{2\gamma-1} \chi_{\{u>0\}} \frac{x - \pi_{\Gamma}(x)}{|x - \pi_{\Gamma}(x)|}\cdot x dx\\ 
& \qquad =  2 \gamma \int_{B_r(0)} \dist(x, \Gamma)^{2\gamma} \chi_{\{u>0\}} dx\\
& \qquad \qquad  + 2 \gamma \int_{B_r(0)} \dist(x, \Gamma)^{2\gamma-1} \chi_{\{u>0\}} \frac{x - \pi_{\Gamma}(x)}{|x - \pi_{\Gamma}(x)|}\cdot \pi_{\Gamma}(x) dx.
\end{align*}
Thus, we obtain the following equation
\begin{align}\label{e:inner variation with eta}\nonumber
0 = & (n + 2\gamma)\int_{B_r(0)}|\nabla u|^2 + Q^{2} \chi_{\{u >0\}}dx  - r\int_{\partial B_r(0)}|\nabla u|^2 + Q^{2} \chi_{\{u >0\}}d\sigma\\
& \qquad - (2 + 2\gamma)\int_{B_r(0)} |\nabla u|^2dx  + 2r\int_{\partial B_r(0)} (\nabla u \cdot \eta)^2 d\sigma\\\nonumber
& \qquad + 2 \gamma \int_{B_r(0)} \dist(x, \Gamma)^{2\gamma-1} \chi_{\{u>0\}} \frac{x - \pi_{\Gamma}(x)}{|x - \pi_{\Gamma}(x)|}\cdot \pi_{\Gamma}(x) dx.
\end{align}

We now using \eqref{e:inner variation with eta} to re-write $U'(r)$.  Recalling \eqref{e: U'}, \eqref{e:inner variation with eta} gives
\begin{align*}
U'(r) = & \frac{1}{r^{n + 1 + 2\gamma}}\biggl(2r \int_{\partial B_r(0)}(\nabla u \cdot \eta)^2d\sigma - (2 + 2\gamma)\int_{B_r(0)}|\nabla u|^2dx\\
& \qquad + 2\gamma \int_{B_r(0)}\dist(x, \Gamma)^{2\gamma-1} \chi_{\{u>0\}} \frac{x - \pi_{\Gamma}(x)}{|x - \pi_{\Gamma}(x)|}\cdot \pi_{\Gamma}(x) dx \biggr).
\end{align*}
Now, using the Divergence theorem one obtains that $\int_{B_r(0)}|\nabla u|^2dx = \int_{\partial B_r(0)}u\nabla u \cdot \eta d\sigma$ for almost every $0<r$.  Whence, for every such radius
\begin{align*}
U'(r) = &\frac{1}{r^{n + 1 + 2\gamma}}\biggl(2r \int_{\partial B_r(0)}(\nabla u \cdot \eta)^2d\sigma - (2 + 2\gamma)\int_{\partial B_r(0)}u \nabla u \cdot \eta d\sigma \\
& \qquad + 2\gamma \int_{B_r(0)}\dist(x, \Gamma)^{2\gamma-1} \chi_{\{u>0\}} \frac{x - \pi_{\Gamma}(x)}{|x - \pi_{\Gamma}(x)|}\cdot \pi_{\Gamma}(x)dx \biggr).
\end{align*}

Now, since, $\frac{d}{dr}W_{\gamma +1}(0, r, u, \Gamma) = U'(r) - (\gamma +1)V'(r),$ we calculate as follows.
\begin{align*}
&\frac{d}{dr}W_{\gamma +1}(0, r, u, \Gamma) \\
& \qquad \qquad  = \frac{1}{r^{n + 1 + 2\gamma}}\biggl(2r \int_{\partial B_r(0)}(\nabla u \cdot \eta)^2d\sigma - (2 + 2\gamma)\int_{\partial B_r(0)}u \nabla u \cdot \eta d\sigma \\
& \qquad \qquad \qquad + 2\gamma \int_{B_r(0)}\dist(x, \Gamma)^{2\gamma-1} \chi_{\{u>0\}} \frac{x - \pi_{\Gamma}(x)}{|x - \pi_{\Gamma}(x)|} \cdot \pi_{\Gamma}(x) dx\biggr)\\
& \qquad \qquad\qquad - (\gamma + 1) \left( -\frac{2 + 2\gamma}{r^{n + 2 +2\gamma}}\int_{\partial B_r(0)}u^2d \sigma + \frac{1}{r^{n +1 +2\gamma}}\int_{\partial B_r(0)}2 u \nabla u \cdot \eta d\sigma\right)\\
& \qquad \qquad = \frac{2}{r^{n + 2\gamma}}\int_{\partial B_r(0)}(\nabla u \cdot \eta - \frac{\gamma + 1}{r}u)^2 \\
& \qquad \qquad \qquad + \frac{1}{r^{n + 1 + 2\gamma}}\left(2\gamma \int_{B_r(0)}\dist(x, \Gamma)^{2\gamma-1} \chi_{\{u>0\}} \frac{x - \pi_{\Gamma}(x)}{|x - \pi_{\Gamma}(x)|}\cdot \pi_{\Gamma}(x) dx\right).
\end{align*}
This proves \eqref{Weiss derivative}.

To finish the lemma, we note that for any $y \in \Gamma$, if $x \in B_r(0)$ is such that $y = \pi_{\Gamma}(x)$, then $x-y \in N_y\Gamma$.  We denote $\frac{x-y}{|x-y|} = \eta_{\Gamma}(\pi_{\Gamma}(x))$.  Hence, we over-estimating as follows.
\begin{align*}
&\left|\frac{2\gamma}{r^{n + 1 + 2\gamma}} \int_{B_r(0)}\dist(x, \Gamma)^{2\gamma-1} \chi_{\{u>0\}} \frac{x - \pi_{\Gamma}(x)}{|x - \pi_{\Gamma}(x)|}\cdot \pi_{\Gamma}(x) dx\right|\\
& \qquad \le \frac{2\gamma}{r^{n + 2}} \int_{B_r(0)} |\eta_{\Gamma}(\pi_\Gamma(x)) \cdot \pi_{\Gamma}(x)| dx\\
& \qquad \le \frac{2\gamma}{r^{n + 2}} \int_{B_r(0)}8[\Gamma]_{\alpha}r^{1 + \alpha} dx \le 16\gamma [\Gamma]_{\alpha}r^{\alpha-1},
\end{align*}
where the penultimate inequality follows from Lemma \ref{l:normal drift}. This gives $\eqref{e:W derivative inequality}$.
\end{proof}

\begin{corollary}\label{c:almost monotonicity}
Let $\Gamma$ be as above, and let $u$ be a local minimizer of $J_{Q}(\cdot, \Omega)$, and $x_0 \in \Gamma \cap \partial \{u>0\}.$  Then, for any $0< r< R< r_0$,
\begin{align}\label{e:almost monotonicity Weiss}
W_{\gamma +1}(x_0, r, u, \Gamma) & \le W_{\gamma +1}(x_0, R, u, \Gamma) + C(\gamma, \alpha) [\Gamma]_{\alpha}R^{\alpha}.
\end{align}
where $c(\gamma, \alpha) = 16\frac{\gamma}{\alpha}.$ Furthermore, $\lim_{r \rightarrow 0^+} W_{\gamma +1}(x_0, r, u, \Gamma) = W_{\gamma +1}(x_0, 0^+, u, \Gamma)$ exists and $W_{\gamma + 1}(x_0, 0^+, u, \Gamma) \in [-c(n,\gamma), c(n)]$.
\end{corollary}

\begin{proof}
The inequality \eqref{e:almost monotonicity Weiss} follows from integrating \eqref{e:W derivative inequality} and rearranging.  To see that the limit exists, we first observe that by Corollary \ref{c:local Lipschitz bound}, 
\begin{align*}
&\frac{1}{r^{n-2+2(\gamma +1)}} \int_{B_r(x_0)}|\nabla u|^2 + Q^2(x)\chi_{\{u>0\}}dx \le r^{-n-2\gamma} \omega_{n}C(n) r^{n+2\gamma} \le C(n) \\
& - \frac{\gamma + 1}{r^{n-1+2\gamma}}\int_{\partial B_r(x_0)} u^2 d\sigma \ge -(\gamma + 1)r^{-n-1-2\gamma} \omega_{n-1} C(\gamma, n)r^{n+1-2\gamma} = -C(n, \gamma).
\end{align*} Now, suppose for the sake of contradiction that $r_i \rightarrow 0$ and $a< b$ are accumulation points of $W_{\gamma+1}(x_0, r_i, u, \Gamma)$.  Then for $r_i$ sufficiently small, $C(\gamma, \alpha)r_i^{\alpha}< |a-b|/2$ and $|W_{\gamma+1}(x_0, r_i, u, \Gamma) - a| < |a-b|/2$ and we obtain a contradiction.
\end{proof}

\begin{remark}\label{r: gamma 0 non 0}(Comparing $\gamma =0$ and $\gamma \not = 0$)
When $\gamma = 0$ and $u$ is a local minimizer of $J_1(\cdot, \Omega)$ and $W_1$ is the Weiss $1$-density associated to $J_{1}(\cdot, \Omega)$, one can verify by direct computation that \begin{align*}
    \frac{d}{dr}W_{1}(x_0, r, u) & = \frac{1}{r^{n}}\int_{\partial B_r(x_0)}\left(\nabla u \cdot \eta - \frac{1}{r} u \right)^2d\sigma + \frac{n}{r^{n+1}}(W_1(x_0, r, z_u) - W_1(x_0, r, u)),
\end{align*} where $z_u$ is the $1$-homogeneous extension of $u|_{\partial B_r(x_0)}$. From here, one typically uses the Noether equations obtained from the inner variation to obtain an ``equipartition of energy" result and rewrites 
\begin{align*}
    \frac{d}{dr}W_{1}(x_0, r, u) & = \frac{2}{r^{n}}\int_{\partial B_r(x_0)}\left(\nabla u \cdot \eta - \frac{1}{r} u \right)^2d\sigma.
\end{align*}
See \cite{Velichkov19} Chapter 9 for details in the $\gamma = 0$ case, which hold \textit{mutatis mutandis} when $0< \gamma$.  

When $0<\gamma$, one may obtain a similar expression, using the $(1+\gamma)$-homogeneous extension. Following the same direct computation and assuming that $\Gamma$ is flat, if $Q(x) = \dist(x, \Gamma)^{\gamma}$, $u$ is a local minimizer of $J_{Q}(\cdot, \Omega)$, and $x_0 \in \Gamma \cap \partial \{u>0\}$ one obtains
\begin{align*}
    & \frac{d}{dr}W_{1+\gamma}(x_0, r, u, \Gamma) = \frac{2}{r^{n +2\gamma}}\int_{\partial B_r(x_0)}\left(\nabla u \cdot \eta - \frac{1+ \gamma}{r} u \right)^2d\sigma\\
    & \qquad + \frac{n+2\gamma}{r^{n+1+ 2\gamma}}(W_{1 + \gamma}(x_0, r, z_u) - W_{1+ \gamma}(x_0, r, u)) + \frac{\gamma(\gamma+1)}{r^{n+1+2\gamma}}\int_{\partial B_r(x_0)}u^2 d\sigma,
\end{align*}
where $z_u$ is now the $(1+ \gamma)$-homogeneous extension of $u|_{\partial B_r(x_0)}$.  Indeed, combined with the weak geometry (Lemma \ref{l:interior balls}) and the Lipschitz bound (Corollary \ref{c:local Lipschitz bound}), this expression is sufficient to obtain the results of Corollary 
\ref{c:almost monotonicity}.  

However, this form is less convenient than \eqref{Weiss derivative} for several reasons.  First, it obscures the role of the geometry of $\Gamma$ which we allow to be non-flat.  While the errors incurred are controlable for non-flat $\Gamma$, for the purposes of this paper it is easier to follow \cite{WeissVarvaruca11} Theorem 3.5, apply the Noether Equations first, and avoid introducing $(1+ \gamma)$-homogeneous extensions. 
\end{remark}

\section{Compactness and Blow-ups}\label{s:companctness}

In this section, we show that $\epsilon_0$-local minimizers of $J_{Q}(\cdot, \Omega)$ satisfy nice compactness properties. In particular, we highlight Lemma \ref{l:positivity set convergence}, which describe the convergence of the positivity sets near $\Gamma = \{Q=0\}$. Note that because of the degeneracy of $Q$, we do not obtain that the positivity sets of convergent subsequences of functions converge in the local Hausdorff metric. See Remark \ref{r:positivity set convergence}.  This behavior is markedly different from known behavior when $\gamma=0$.  The rest of the results follow from standard techniques.

As a consequence of the compactness result (Theorem \ref{t:compactness}, below) the second collection of results in this section are results on the existence and characterization of blow-ups.  In particular, in Lemma \ref{l:homogeneity} it is shown that if $Q(x) = \dist(x, \Gamma)^{\gamma}$ then blow-ups of local minimizers at points in $\Gamma \cap \partial \{ u>0\}$ are $(1+\gamma)$-homogeneous.

First, we define the following rescalings.

\begin{definition}(Rescalings)
Let $n \ge 2$.  For any function $f: \mathbb{R}^n \rightarrow \mathbb{R}$, $x \in \mathbb{R}^n$, and $0< r< \infty$ we define the rescalings
\begin{align*}
f_{x, r}(y) := \frac{f(ry + x)}{r^{\gamma +1}}.
\end{align*}
For any set $K\subset \mathbb{R}^n$ we define the rescalings, 
\begin{align*}
K^{x, r}: = \frac{1}{r}(K - x)
\end{align*}
\end{definition}

\begin{remark}
We define these rescalings in addition to the rescalings of \eqref{def T rescalings}, because whereas the rescalings $T_{x, r}u$ are used to define the quantitative stratification, they do not obviously work with the structure of local minimizers.  If $0 \le k\le n-1$ is an integer, $\Gamma$ is a $k$-dimensional $C^{1, \alpha}$ submanifold, $0< \gamma$, and $Q(x) = \dist(x, \Gamma)^{\gamma}$ then by change of variables it is clear that if $u$ is an $\epsilon_0$-minimizer of $J_{Q}(\cdot, \Omega)$, then 
$u_{x, r}$ is an $\epsilon_0 \min\{r^{-2\gamma}, r^{-n}\}$-minimizer of $J_{Q_{x, r}}(\cdot, \Omega^{x, r})$.  Furthermore, if $Q^{x, r}(y) = \dist(y, \Gamma^{x, r})^{\gamma}$.
\end{remark}

\begin{remark}\label{r:manifold convergence}
We note that since $\Gamma$ is a $k$-dimensional $C^{1, \alpha}$ submanifold, for any sequence $r_j \rightarrow r \in [0, \infty)$ and any $x_0 \in \Gamma$, the sequences $\Gamma^{x_0, r_j} \rightarrow \Gamma^{x_0, r}$ locally in the Hausdorff metric on compact subsets.  When $r =0$, we will denote $\Gamma^{x_0, 0}$ by $\Gamma^{\infty}$ and note that $\Gamma^{\infty} = T_{x_0}\Gamma$ a $k$-dimensional linear subspace of $\mathbb{R}^n$.
\end{remark}

\begin{theorem}\label{t:compactness}(Compactness)
Let $n\ge 2$ and $0 \le k \le n-1$ be integers. Let $0< \epsilon_0$ and $0<\gamma<\infty$.  Let $\Gamma_i$ be a $(1, M)$-$C^{1, \alpha}$ submanifold satisfying $0 \in \Gamma_i$.  Let $Q_i(x) = \dist(x, \Gamma_i)^{\gamma}$, and let $\Omega_i$ be a Lipschitz domain with $\overline{B_2(0)} \subset \Omega_i$.

Suppose that $\{u^i\}_i$ is a sequence of $\epsilon_0$-local minimizer of $J_{Q_i}(\cdot, \Omega_i)$ with standard scale $r_{0, i}(\epsilon_0, u_i, \Omega_i) = 1$ which satisfies 
$$
0 \in \partial \{u_i>0\} \cap \Gamma_i.
$$
Furthermore, let $\{r_i \}_i$ be a sequence of numbers $0<r_i \le 2$ such that and $r_i \rightarrow r \in [0, 2)$.  Write $B= B_\frac{1}{r}(0)$ if $r>0$ and $B = \mathbb{R}^n$ if $r = 0$. 

Then, there is a subsequence, $r_j \rightarrow r$, a $k$-dimensional $(1, M)$-$C^{1, \alpha}$ submanifold $\Gamma$, and a function $u \in C^{0, 1}_{loc}(B) \cap W^{1, 2}_{loc}(B)$ such that,
\begin{enumerate}
\item $u^j_{0, r_j} \rightarrow u$ in $C^{0, 1}_{loc}(B) \cap W^{1, 2}_{loc}(B)$.
\item For every $0< R< \diam(B)$, 
\begin{align*}
\Gamma_j^{0, r_j} \cap B_R(0) \rightarrow \Gamma \cap B_R(0)
\end{align*}
in the Hausdorff metric on compact subsets.
\item For every $0<R< \diam(B)$ and any $\epsilon >0$,  there is an $N \in \mathbb{N}$ such that for all $j \ge N$,
\begin{align*}
\partial \{u > 0 \} \cap B_R(0) \subset B_{\epsilon}(\partial \{u^j_{0, r_j} > 0 \}).
\end{align*}
\item For any $0< R< \diam(B)$ and any $\epsilon > 0$, there is an $N \in \mathbb{N}$ such that for all $j \ge N$,  
\begin{align*}
\partial \{u^j_{0, r_j} > 0 \} \cap B_R(0) \cap \{u > 0 \} \subset B_{\epsilon}(\partial \{u > 0 \}).
\end{align*}
Similarly, for any $0< R< \diam(B)$ and any $\epsilon > 0$, there is an $N \in \mathbb{N}$ such that for all $j \ge N$,  
\begin{align*}
\partial \{u^j_{0, r_j} > 0 \} \cap B_R(0) \cap \{u = 0 \} \subset B_{\epsilon}(\partial \{u > 0 \} \cup \Gamma).
\end{align*}
\item $\chi_{\{u^j_{0, r_j} >0 \}} \rightarrow \chi_{\{u>0 \}}$ in $L^1_{loc}(B)$.
\item The function $u$ is harmonic in $\{u > 0\}$.
\end{enumerate}
\end{theorem}

The proof of Theorem \ref{t:compactness} is broken up into the following lemmata.  We note that $(2)$ is given as Lemma \ref{l:manifold convergence}, and that $(5)$ follows immediately from $(3)$ and $(4)$ and the fact that $\Gamma$ is a lower-dimensional submanifold.

\begin{lemma}\label{l:Rellich compactness}
Let $u^i$, $\Gamma_i$, $r_i \rightarrow r$, and $B$ as in Theorem \ref{t:compactness}. There is a subsequence $r_j \rightarrow r$ and a function $u \in C^{0, 1}_{loc}(B) \cap W^{1, 2}_{loc}(B)$ such that
\begin{enumerate}
\item $u^j_{0, r_j} \rightarrow u$  in  $L^{2}_{loc}(B)$ and $C^{0, 1}_{loc}(B)$.
\item $\nabla u^j_{0, r_j} \rightharpoonup \nabla u$ in $L^2_{loc}(B ; \mathbb{R}^n).$
\end{enumerate}
\end{lemma}

\begin{proof}  Fix a $0<R< \diam(B).$  Since $0 \in \Gamma_i \cap \partial \{u^i>0\}$, the local Lipschitz bounds in Corollary \ref{c:local Lipschitz bound} immediately implies that,
\begin{align*}
||\nabla u^i_{0, r_i}||_{L^2(B_R(0))} & = \left( \int_{B_R(0)}\left(\frac{\nabla u^i(r_iy)}{r_i^{\gamma}}\right)^2dy \right)^{\frac{1}{2}} \le C(n)R^{\gamma} R^{\frac{n}{2}}.
\end{align*}
Thus, $\nabla u^i_{0, r_i}$ are uniformly bounded in $L^2(B_R(0); \mathbb{R}^n)$. It remains to show that $u^i_{0, r_i}$ are uniformly bounded in $L^2(B_R(0))$. By Corollary \ref{c:local Lipschitz bound} and the assumption that $0 \in \Gamma_i \cap \partial \{u^i >0\},$ for any $0< r$ we bound
\begin{align*}
\int_{B_r(0)}|u^i|^2dx \le C(n)r^{n+2(\gamma + 1)}.
\end{align*}
Thus, $||u^i_{0, r_i}||_{L^2(B_R(0))} \le CR^{\gamma +1}R^{\frac{n}{2}}.$  By Rellich compactness, then, there exists a function $u$ and a subsequence $r_j \rightarrow 0$ such that, 
\begin{align*}
u^j_{0, r_j} \rightarrow u & \text{  in } L^{2}_{loc}(B)\\
\nabla u^j_{0, r_j} \rightharpoonup \nabla u & \text{  in } L^{2}_{loc}(B;\mathbb{R}^n).
\end{align*}
Note that by the uniform local Lipschitz bounds of $u^j_{0, r_j},$ we have that $u^j_{0, r_j}$ is locally Lipschitz continuous.  Thus, up to a further subsequence $u^j_{0, r_j} \rightarrow u$ in $C^{0, 1}(B_R(0))$ by Arzela-Ascoli.  Diagonalizing for a sequence of $R \nearrow \diam(B)$ proves the full claim.
\end{proof}

\begin{lemma}\label{l:positivity set convergence}
For $u^i$, $\Gamma_i$, $r_i \rightarrow r$, and $B$ as in Theorem \ref{t:compactness}, the following holds.  Let $u^j_{0,r_j}$ be a subsequence as in Lemma \ref{l:Rellich compactness} such that $\Gamma_j^{0, r_j} \rightarrow \Gamma$ as in Theorem \ref{t:compactness}(2).  

For every $0<R< \diam(B)$ and any $\epsilon >0$,  there is an $N \in \mathbb{N}$ such that for all $j \ge N$, such that,
\begin{align*}
\partial \{u > 0 \} \cap B_R(0) & \subset B_{\epsilon}(\partial \{u^j_{0, r_j} > 0 \})\\
\partial \{u^j_{0, r_j} > 0 \} \cap B_R(0) \cap \{u > 0 \}& \subset B_{\epsilon}(\partial \{u > 0 \})\\
\partial \{u^j_{0, r_j} > 0 \} \cap B_R(0) \cap \{u = 0 \} & \subset B_{\epsilon}(\partial \{u > 0 \} \cup \Gamma).
\end{align*}
\end{lemma}

\begin{remark}\label{r:positivity set convergence}
We note that we do not have convergence of the free boundaries locally in the Hausdorff metric on compact subsets.  It is \textit{a priori} possible that the positivity sets, $\{u^j_{0, r_j}>0 \}$ have tenderils which reach out into $\{u = 0 \}$ and cleave close to $\Gamma_j$ which thin as $j \rightarrow \infty$ and vanish in the limit.  This problem is connected to the degeneracy of $Q$ and the problem of establishing interior balls in balls centered on $\Gamma \cap \partial \{u_{0}>0 \}.$
\end{remark}

\begin{proof}
Let $0< R < \diam(B)$ and $1\ge \epsilon >0$ be given.  We prove the first containment.  Assume for the sake of contradiction that $x' \in \partial \{u>0 \},$ but that $x' \not \in B_{\epsilon}(\partial \{u^j_{0, r_j}>0 \})$ for all $j \in \mathbb{N}$.  Suppose there is a subsequence $j'$ such that $x' \in \{u^{j'}_{0, j'}>0\}$. By Corollary \ref{c:u growth}, we have that $u^j_{0, r_j}(x') \ge C(n, \gamma) \epsilon^{1 + \gamma}.$

Since $x' \in \partial \{u>0 \}$ and $u^j_{0, r_j} \rightarrow u$ in $C^{0, 1}_{loc}(B),$ we have that 
\begin{align*}
 0 = u(x') = \lim_{j' \rightarrow \infty}u^{j'}_{0, j'}(x') \ge C(n, \gamma) \epsilon^{1+\gamma} > 0.
\end{align*}
This is absurd.  On the other hand, there exists a subsequence $j'$ such that $x' \in \{u^{j'}_{x_0, j'} = 0 \},$ we obtain a similar contradiction.  Indeed, since $x' \in \partial \{u>0 \}$ there must exist some $y \in B_{\epsilon}(x')$ such that $u(y)>0$. Therefore, by $C^{0, 1}$ convergence, there must exist a number $N$ sufficiently large such that $u^{j'}_{0, j'}(y)>0$ for all $j' \ge N$. This forces $B_{\epsilon}(x') \cap \partial \{u^j_{0, r_j}>0 \} \not = \emptyset$ and proves the claim.

To see the second containment result, we assume for the sake of contradiction that there exists a sequence of $j \in \mathbb{N}$, such that there exists an $x_j \in B_R(0) \cap \partial \{u^j_{0, r_j}> 0\} \cap \{u> 0\}$ for which
\begin{align*}
x_j \not \in B_{\epsilon}(\partial \{u> 0 \}).
\end{align*}  
By passing to a further subsequence, we may assume that $x_j \rightarrow x' \in \overline{B_R(0)} \cap \{u>0 \}.$  Since $u^j_{0, r_j} \rightarrow u$ in $C^{0, 1}_{loc}(\mathbb{R}^n)$ which forces $0< u(x') = \lim_{j \rightarrow \infty}u^j_{0, r_j}(x_j) = 0.$ This is a contradiction. Hence, the claim follows.

The third containment follows from an analogous argument.  We assume for the sake of a contradiction that exists a sequence of $j \in \mathbb{N}$, for which we can find a point $x_j \in B_R(0) \cap \partial \{u^j_{0, r_j}> 0\} \cap \{u= 0\}$ for which
\begin{align*}
x_j \not \in B_{\epsilon}(\partial \{u> 0 \} \cup \Gamma^{0, r_j}_j).
\end{align*}  
However, by Lemma \ref{l:interior balls}, and our choice of $Q$, we have that for all $x_j$, there is a ball $B_{c\epsilon}(y_j) \subset B_{\epsilon/2}(x_j) \cap \{u^j_{0, r_j}>0 \}.$  Passing to a further subsequence $j'$, we may assume that $x_{j'} \rightarrow x'$ and $y_{j'} \rightarrow y'$.  By $C^{0, 1}$ convergence, then, $B_{c\epsilon}(y') \subset \{u>0 \}.$  This contradicts the assumptions that $x_{j'} \in \{u=0\}$ and $\dist(x_{j'}, \partial \{u >0 \}) > \epsilon$ for all $j'$.
\end{proof}

\begin{lemma}\label{l:harmonic limit}Let $u^i,$ $r_j$, and $u$ be as in Theorem \ref{t:compactness}.  Then $u$ is harmonic in $\{u >0\}.$   
\end{lemma}

\begin{proof}
Let $x' \in \{u>0 \}$ with $\dist(x', \partial \{u>0 \}) = 3\epsilon$. For sufficiently large $j$, we have that $B_{\epsilon}(x') \subset \{u_{0, r_j}> 0 \}.$  Since $u^j_{0, r_j}$ is harmonic in $B_{\epsilon}(x'),$ and the $W^{1, 2}(B_{2\epsilon}(x'))$ norm of $u^j_{0, r_j}$ are uniformly bounded, it is a classical result that $u^j_{0, r_j}$ converge in $C^{\infty}(B_{2\epsilon}(x'))$.  Thus, the harmonicity of $u$ at $x'$ is a consequence of $u$ satisfying the mean value property in $B_{\epsilon}(x').$  That is, for all $y' \in B_{\epsilon}(x'),$
\begin{align*}
u(y') = \lim_{j \rightarrow \infty} u^j_{0, r_j}(y') = \lim_{j \rightarrow \infty} \dashint_{B_{\epsilon}(y')}u^j_{0, r_j}d\sigma = \dashint_{B_{\epsilon}(y')}u d\sigma.
\end{align*}
\end{proof}

\begin{lemma}\label{l:stong convergence}
Let $u^i$, $u$, and $r_m$ be as in Theorem \ref{t:compactness}.  Then there exists a subsequence such that $u^j_{0, r_j} \rightarrow u$ strongly in $W^{1, 2}_{\loc}(\mathbb{R}^n).$
\end{lemma}

\begin{proof}
By Lemma \ref{l:Rellich compactness}, we may reduce to a subsequence such that $u^j_{x_0, r_j} \rightharpoonup u$ in $W^{1, 2}_{\loc}(\mathbb{R}^n).$  To show strong convergence, it suffices to show that for every $\phi \in C^{\infty}_{c}(\mathbb{R}^n)$
\begin{align*}
\limsup_{j \rightarrow \infty} \int |\nabla u^j_{x_0, r_j}|^2 \phi dx \le \int |\nabla u|^2\phi dx.
\end{align*}
Using the fact that $u$ and $u^j_{x_0, r_j}$ are harmonic in their positivity sets and the senses of convergence of $u^j_{x_0, r_j} \rightarrow u$ in Lemma \ref{l:Rellich compactness}, we obtain by integration by parts,
\begin{align*}
\int |\nabla u^j_{x_0, r_j}|^2\phi dx & = -\int u^j_{x_0, r_j}\nabla u^j_{x_0, r_j}\cdot \nabla \phi dx\\
& \rightarrow -\int u \nabla u\cdot \nabla \phi dx = \int |\nabla u|^2\phi dx.
\end{align*}
It follows that $u^j_{x_0, r_j} \rightarrow u$ strongly in $W^{1, 2}_{\loc}(\mathbb{R}^n).$
\end{proof}

This concludes the proof of Theorem \ref{t:compactness}.

\begin{definition}\label{d:blow-ups}(Blow-ups)
Let $n\ge 2$ and $0 \le k \le n-1$ be integers. Let $0<\gamma<\infty$.  Let $\Gamma$ be a $(1, M)$-$C^{1, \alpha}$ submanifold satisfying $0 \in \Gamma$.  Let $Q(x) = \dist(x, \Gamma)^{\gamma}$, and let $\Omega$ be a Lipschitz domain.

If $u$ is a local minimizer of $J_{Q}(\cdot, \Omega)$, $x_0 \in \Gamma \cap \partial\{u> 0\} \cap \Omega$, and $r_i \rightarrow 0$ any subsequential limit $u_{x_0, r_j} \rightarrow u_{x_0, 0}$ in the senses of Theorem \ref{t:compactness} is called a \textit{blow-up} of $u$ at $x_0$.   
\end{definition}

\begin{corollary}\label{c:Weiss limit}
Let $u$ be a local minimizer of $J_{Q}(\cdot, B_1(0))$ for $\Gamma$ a $k$-dimensional $(1, M)$-$C^{1, \alpha}$ submanifold with $x_0 \in \Gamma \cap \partial \{u>0\} \cap B_1(0)$.  Let $r_j \rightarrow 0$ and $u_{x_0, r_j} \rightarrow u_{x_0, 0}$ in the senses of Theorem \ref{t:compactness}.  Then for any fixed $0<r< \infty$,
\begin{align*}
\lim_{j \rightarrow \infty}W_{\gamma +1}(0, r, u_{x_0, r_j}, \Gamma^{x_0, r_j}) = W_{\gamma +1}(0, r, u_{x_0, 0}, T_{x_0}\Gamma).
\end{align*}
\end{corollary}

The corollary follows immediately from Theorem \ref{t:compactness}(1) and (5).

\begin{lemma}(Homogeneity) \label{l:homogeneity}
Let $x_0 \in \mathcal{S}$ for a local minimizer $u$, and let $u_{x_0, 0}$ be a blow-up as in Definition \ref{d:blow-ups}. The function $u_{x_0, 0}$ is $(\gamma +1)$-homogeneous and non-trivial.
\end{lemma}

\begin{proof}
We recall that the Weiss density is invariant under the dilation $W_{\gamma+1}(x_0, r, u, \Gamma) = W_{\gamma + 1}(0, 1, u_{x_0, r}, \Gamma^{x_0, r})$.  Therefore, for all $0< r_1 < r_2 < \infty$ and all sufficiently small $0< r< \dist(0, \partial \Omega)$,
\begin{align*}
& W_{\gamma + 1}(0, r_1, u_{x_0, r_j}, \Gamma^{0, r_j}) - W_{\gamma + 1}(0, r_2, u_{x_0, r_j}, \Gamma^{0, r_j}) = \\
& \qquad \int_{r_1}^{r_2} \frac{1}{s^{n + 2+ 2\gamma}}\int_{\partial B_s}(\nabla u_{0, r_j}(x)\cdot x - (\gamma +1)u_{0, r_j}(x))^2d\sigma(x) ds\\
& \qquad + \int_{r_1}^{r_2}\frac{1}{s^{n + 1 + 2\gamma}}2\gamma \int_{B_s(0)}\dist(x, \Gamma^{0, r_j})^{2\gamma-1} \chi_{\{u_{0, r_j}>0\}} \frac{(x - \pi_{\Gamma^{0, r_j}}(x))}{|x - \pi_{\Gamma^{0, r_j}}(x)|}\cdot \pi_{\Gamma^{0, r_j}}(x) dxds.
\end{align*}
We now calculate the limit as $j \rightarrow \infty$, in two ways.  First, we recall that $\Gamma^{0, r_j}$ converges to a $k$-dimensional linear subspace locally in the Hausdorff metric.  Thus, the function 
\begin{align*}
\frac{(x - \pi_{\Gamma^{0, r_j}}(x))\cdot \pi_{\Gamma^{0, r_j}}(x)}{|x - \pi_{\Gamma^{0, r_j}}(x)|} \rightarrow 0
\end{align*}
in $L^1(B_1(0))$.  Recalling Corollary \ref{c:Weiss limit} and taking the limit as $j \rightarrow \infty$, then, we obtain the following expression.
\begin{align*}
& W_{\gamma + 1}(0, r_1, u_{x_0, 0}, \Gamma^{x_0, 0}) - W_{\gamma + 1}(0, r_2, u_{x_0, 0}, \Gamma^{x_0, 0})  \\
& \qquad = \int_{r_1}^{r_2} \frac{1}{s^{n + 2+ 2\gamma}}\int_{\partial B_s}(\nabla u(x)\cdot x - (\gamma +1)u(x))^2d\sigma(x) ds.
\end{align*}
We now show that this expression is zero. Because for any sequence $r_j \rightarrow 0$,
\begin{align*}
    \lim_{j \rightarrow \infty} W_{\gamma +1}(x_0, r_j, u, \Gamma) = W_{\gamma +1}(x_0, 0^+, u, \Gamma),
\end{align*}
it must be the case that for any $c \in (0,1]$,
\begin{align*}
\lim_{j \rightarrow \infty} W_{\gamma +1}(x_0, c r_j, u, \Gamma) & =     \lim_{j \rightarrow \infty} W_{\gamma +1}(0, c, u_{x_0, r_j}, \Gamma^{x_0, r_j})\\
& = W_{\gamma +1}(0, c, u_{x_0, 0}, T_{x_0}\Gamma)
\end{align*}
where the last equality is obtained by Corollary \ref{c:Weiss limit}. Comparing the two limits proves that 
\begin{align*}
0 &= W_{\gamma + 1}(0, r_1, u_{x_0, 0}, \Gamma^{x_0, 0}) - W_{\gamma + 1}(0, r_2, u_{x_0, 0}, \Gamma^{x_0, 0})  \\
& \qquad = \int_{r_1}^{r_2} \frac{1}{s^{n + 2+ 2\gamma}}\int_{\partial B_s}(\nabla u_{x_0, 0}(x)\cdot x - (\gamma +1)u_{x_0, 0}(x))^2d\sigma(x) ds.
\end{align*}
Thus, $u_{x_0, 0}$ is $(\gamma + 1)$ homogeneous in $B_{r_2}(0) \setminus B_{r_1}(0).$  Repeating this for $r_1 \rightarrow 0$ and $r_2 \rightarrow \infty$ shows that $u_{x_0, 0}$ is $(\gamma + 1)$ homogeneous in $\mathbb{R}^n$. It is non-trivial because $x_0 \in \mathcal{S}$ and so by assumption $\{u_{x_0, 0}>0\}$ is non-trivial.
\end{proof}

\begin{corollary}\label{c:limit solution}
Fix $0< \gamma$ and $0 \le k \le n-1$.  Let $\overline{B_2(0)} \subset \Omega _i$, and let $u_i$ be a sequence of $\epsilon_0$-local minimizers of $J_{Q_i}(\cdot, \Omega_i)$ such that $Q_i(x) = \dist(x, \Gamma_i)^\gamma$ for $\Gamma_i$ a $k$-dimensional $(1, \frac{1}{4})$-$C^{1, \alpha}$ submanifold with $0 \in \Gamma_i$. Assume that the standard scales $r_{0,i} = 1$.  For a subsequence $u_i \rightarrow u$ and their corresponding $\Gamma_i \rightarrow \Gamma$ in the senses of Theorem \ref{t:compactness}, the limiting function $u$ satisfies the Noether equations (\ref{e:Noether equations}) and the almost monotonicity equations of Lemma \ref{t:Weiss almost monotonicity} and Corollary \ref{c:almost monotonicity}.
\end{corollary}

\begin{proof}
By assumption, $u_i \rightarrow u$ strongly in $W^{1, 2}(B_1(0)) \cap C^{0, 1}(B_1(0))$ and $\Gamma_i \cap B_1(0) \rightarrow \Gamma \cap B_1(0)$ in the Hausdorff metric on compact subsets.  The first point implies that $|\nabla u_i|^2 \rightarrow |\nabla u|^2$ in $L^{1}(B_1(0))$.  The second point implies that $\dist(\cdot, \Gamma_i) \rightarrow \dist(\cdot,\Gamma)$ in $L^{1}(B_1(0)) \cap C^{0, 1}(B_1(0)).$  Therefore, $\dist(\cdot, \Gamma_i)^{2\gamma} \rightarrow \dist(\cdot, \Gamma)^{2\gamma}$ in $L^1(B_1(0))$.  Thus, to show that $u$ satisfies the Noether Equations (\ref{e:Noether equations}) it is sufficient to show that $\dist(\cdot, \Gamma_i)^{2\gamma -1} \rightarrow \dist(\cdot, \Gamma)^{2\gamma -1}$ in $L^{1}(B_1(0))$ and $\nabla(\dist(\cdot, \Gamma_i)) \rightharpoonup \nabla(\dist(\cdot, \Gamma))$ in $L^1(B_1(0)).$

We first we note that for any $k$-dimensional $(1, \frac{1}{4})$-$C^{1, \alpha}$ submanifold $\Gamma_i$ satisfying $0 \in \Gamma$, using Lemma \ref{l:distance estimate}, we may estimate
\begin{align*}
\int_{B_1(0)}\dist(x, \Gamma_i)^{2\gamma -1}dx & \le C(k, n)\sum_{i = 1}^{\infty}2^{-i(n-k)}2^{(-i+1)(2\gamma -1)} < \infty\\
\int_{B_{2^{-N}}(\Gamma_i)}\dist(x, \Gamma_i)^{2\gamma -1}dx & \le C(k, n)\sum_{i = N}^{\infty}2^{-i(n-k)}2^{(-i+1)(2\gamma -1)}
\end{align*}
Thus, for any $k$-dimensional $(1, \frac{1}{4})$-$C^{1, \alpha}$ submanifold $\Gamma_i$ satisfying $0 \in \Gamma_i$, given any $\epsilon >0$, there is $\delta> 0$ such that $\int_{B_{\delta}(\Gamma_i)}\dist(x, \Gamma)^{2\gamma -1}dx < \epsilon,$ and that this $\delta= \delta(\epsilon)$.  That is, this estimate is independent of $i \in \mathbb{N}.$  

Now, let $\epsilon > 0$ and let $\delta = \delta(\epsilon)>0$.  Since $\Gamma_i \rightarrow \Gamma$ locally in the Hausdorff metric on compact subsets, for all sufficiently large $i \in \mathbb{N}$, we have that $B_{\frac{1}{2}\delta}(\Gamma_i) \subset B_{\delta}(\Gamma)$.  In particular, for all sufficiently large $i \in \mathbb{N}$, both functions $\dist(\cdot, \Gamma_i)^{2\gamma -1}$ and $\dist(\cdot, \Gamma)^{2\gamma -1}$ are uniformly Lipschitz in $B_1(0) \setminus B_{\delta}(\Gamma)$ and converge in $L^1(B_1(0)).$ Hence,
\begin{align*}
&\lim_{i \rightarrow \infty}||\dist(\cdot, \Gamma_i)^{2\gamma -1}- \dist(\cdot, \Gamma)^{2\gamma -1}||_{L^1(B_1(0))}\\
& \qquad \qquad = \lim_{i \rightarrow \infty}\int_{B_1(0) \cap B_{\delta}(\Gamma)}|\dist(x, \Gamma_i)^{2\gamma -1}- \dist(x, \Gamma)^{2\gamma -1}|dx \\
& \qquad \qquad \qquad + \lim_{i \rightarrow \infty}\int_{B_1(0) \setminus B_{\delta}(\Gamma)}|\dist(x, \Gamma_i)^{2\gamma -1}- \dist(x, \Gamma)^{2\gamma -1}|dx\\
& \qquad \qquad \le 2\epsilon.
\end{align*}
Since $\epsilon>0$ was arbitrary, $\dist(\cdot, \Gamma_i)^{2\gamma -1} \rightarrow \dist(\cdot, \Gamma)^{2\gamma -1}$ in $L^{1}(B_1(0))$.

That $\nabla(\dist(\cdot, \Gamma_i)) \rightharpoonup \nabla(\dist(\cdot, \Gamma))$ in $L^1(B_2(0);\mathbb{R}^n)$ follows immediately from Rellich-Kondrakov compactness, noting that $\nabla(\dist(\cdot, \Gamma_i))$ are uniformly bounded in $W^{1, 2}(B_2(0)).$ This concludes the proof that $u$ satisfies the Noether equations (\ref{e:Noether equations}).  Following the proof of Lemma \ref{t:Weiss almost monotonicity} and Lemma \ref{c:almost monotonicity} verifies that $u$ satisfies the almost monotonicity equations \eqref{Weiss derivative}, \eqref{e:W derivative inequality}, and \eqref{e:almost monotonicity Weiss} with $[\Gamma]_{\alpha} \le \limsup_{i \rightarrow \infty} [\Gamma_i]_{\alpha} \le \frac{1}{4}.$
 \end{proof}

\section{Density lower bounds}\label{s:density}
We now prove some important results about the Weiss density.  One of the key results in an energy ``gap" result in Lemma \ref{l:energy lower-bound}. This next lemma justifies the notion of $W_{\gamma +1}$ as a ``density."

\begin{lemma}\label{l:density}
Fix $0< \gamma$ and $0 \le k \le n-1$.  Let $Q(x) = \dist(x, \Gamma)^\gamma$ for $\Gamma$ a $k$-dimensional $(1, \frac{1}{4})$-$C^{1, \alpha}$ submanifold with $0 \in \Gamma$.  Let $u$ be a local minimizer of $J_{Q}(\cdot, \Omega)$.  Let $x_0 \in \Gamma \cap \partial \{u>0 \}$.  Then,
\begin{align*}
W_{\gamma +1}(x_0, 0^+, u, \Gamma) = \lim_{r \rightarrow 0^+}\frac{1}{r^{n + 2\gamma}}\int_{B_r(x_0)}Q^{2}(x)\chi_{\{u>0 \}}dx.
\end{align*}
Therefore, $W_{\gamma +1}(x_0, 0^+, u, \Gamma) \in [0, c_n],$ where $c_n = \int_{B_1(0)}\dist(x, T_{x_0}\Gamma)^{2\gamma}dx.$ In particular, $W_{\gamma +1}(x_0, 0^+, u, \Gamma) =0$ implies that every blow-up $u_{x_0, 0}\equiv 0.$ Furthermore, the function $x \mapsto W_{\gamma +1}(x, 0, u, \Gamma)$ is upper semicontinuous when restricted to $\Gamma \cap \partial \{u>0 \}.$
\end{lemma}

\begin{proof}
By translation and scaling, we may assume that $x_0 = 0$ and that the standard scale of $\epsilon_0$-local minimizers in $K_{u, \Omega} = 1$. Therefore, Theorem \ref{t:compactness} implies there exists a function $v$ and a sequence of radii $r_j \rightarrow 0$, such that $u_{x_0, r_j}\rightarrow v$ strongly in $W^{1, 2}_{loc}(\mathbb{R}^n)$ and in $C^{0, 1}_{loc}(\mathbb{R}^n)$.  Thus, we calculate,
\begin{align*}
& \lim_{r \rightarrow 0^+}W_{\gamma +1}(0, 1, u_{x_0, r}, \Gamma^{x_0, r}) = \lim_{j \rightarrow \infty}W_{\gamma +1}(0, 1, u_{x_0, r_j}, \Gamma^{x_0, r_j}) \\
&\qquad \qquad  = \int_{B_1(0)}|\nabla v|^2 - (\gamma + 1)\int_{\partial B_1(0)}v^2 d\sigma + \lim_{j \rightarrow \infty} \int_{B_1(0)}\dist(x, \Gamma^{x_0, r_j})^{2\gamma}\chi_{\{u_{x_0, r_j}>0 \}}dx\\
& \qquad \qquad =  \lim_{j \rightarrow \infty} \frac{1}{r^{n + 2\gamma}}\int_{B_{r_j}(x_0)}\dist(x,\Gamma)^{2\gamma}\chi_{\{u>0 \}}dx.
\end{align*}
where the last inequality follows from the fact that $v$ is a $(\gamma + 1)$-homogeneous function. We note that while the limit function $v$ may \textit{a priori} depend upon the subsequence $r_j,$ Lemma \ref{t:Weiss almost monotonicity} implies that the limit on the left-hand side exists.  Therefore, the limit on the right-hand side is unique, even if $v$ may not be.

Upper semicontinuity, in general, holds for limits of monotone increasing functions.  We prove it here for the almost monotone Weiss density. Let $x_0 \in \partial \{u>0 \} \cap B_1(0).$  Let $\delta >0$ and let $0< r(x, \delta)$ be such that $W_{\gamma +1}(x_0, r, u, \Gamma) \le W_{\gamma +1}(x_0, 0, u, \Gamma) + \delta.$  Then for $x \in \Gamma \cap \partial \{u>0 \} \cap B_1(0)$ such that $\mathcal{H}^n(B_r(x) \Delta B_r(x_0)) \le \delta r^{n}$, $|x - x_0| \le \delta$, and $r \le \delta^{\frac{1}{\alpha}},$
\begin{align*}
W_{\gamma+1}(x, 0^+, u, \Gamma) & \le W_{\gamma +1}(x, r, u, \Gamma) + C[\Gamma]_{\alpha}r^{\alpha}\\
 & \le W_{\gamma +1}(x_0, r, u, \Gamma) + |B_r(x) \Delta B_r(x_0)|\frac{Cr^{2\gamma}}{r^{n-2+2(\gamma+1)}}\\
& \qquad \qquad + C\frac{\omega_{n-1}r^{n-1}}{r^{n-1+2(\gamma +1)}}|x-x_0|^{2(\gamma+1)} + C[\Gamma]_{\alpha}r^{\alpha}\\
 & \le W_{\gamma +1}(x_0, r, u, \Gamma) + C\delta\\
& \le W_{\gamma +1}(x_0, 0, u, \Gamma) + (C + 1)\delta.
\end{align*}
 Thus, $\limsup_{x \rightarrow x_0}W_{\gamma +1}(x, 0, u, \Gamma) \le W_{\gamma + 1}(x_0, 0, u, \Gamma).$
\end{proof}

\begin{lemma}\label{l:energy lower-bound}(Density lower bound)
Let $\Gamma$ be a $k$-dimensional $(1, \frac{1}{4})$-$C^{1, \alpha}$ submanifold.  Let $u$ be a local minimizer of $J_{Q}(\cdot, B_2(0))$ for $Q(x) = \dist(x, \Gamma)^{\gamma}$.  Let $x \in \partial \{ u>0\} \cap \Gamma$.  There is a constant, $0<c_0(n, \gamma)$, such that if $W_{\gamma+1}(x, 0^+, u, \Gamma) \not =0$, then 
\begin{align*}
W_{\gamma +1}(x, 0^+, u, \Gamma) \ge c_0.
\end{align*}
\end{lemma}

\begin{proof}
The proof relies upon non-degeneracy of the functions $u$, and in particular the iterior ball condition of Lemma \ref{l:interior balls}.  We note that by Lemma \ref{l:homogeneity}, every blow-up $u_{x, 0^+}$ is homogeneous, and $\{u_{x, 0^+}>0 \}$ is a cone.  Therefore, we observe that if $\{u_{x, 0^+}>0 \} \cap B_1(0) \not \subset B_{c}(T_x\Gamma)$ for some $c = c(n, \gamma)>0$, then Lemma \ref{l:interior balls} and Lemma \ref{l:density} gives the desired result. The strategy of the proof will be to show that such a constant $c(n, \gamma)>0$ must exist.

We argue by contradiction. Let $L$ be any $(n-1)$-dimensional linear subspace which contains $T_x\Gamma$. Suppose that $\{u_{x, 0^+}>0 \} \cap B_1(0)$ is compactly contained within the tubular neighborhood $B_{c}(T_x\Gamma) \subset B_c(L)$.  Define the region $D_c$ as follows,
\begin{align*}
    D_c = B_1(0) \cap \{(r, \theta) \in \mathbb{R}_+ \times \mathbb{S}^{n-1}: (1, \theta) \in \partial B_1(0) \cap B_c(L) \}.
\end{align*}
That is, $D_c$ is the intersection of the cone over $\partial B_1(0) \cap B_c(L)$ and $B_1(0).$ We note that $\{u_{x, 0^+}>0 \} \cap B_1(0) \subset D_c.$  

We define an auxiliary function $h$ as the solution to the following Dirichlet problem, 
\begin{align*}
    \begin{cases}
    \Delta h = 0 & \text{  in  } D_c\\
    h = 0 & \text{  on  } \partial D_c \cap B_1(0)\\
    h = C & \text{  on  } \partial D_c \cap \partial B_1(0)
    \end{cases}
\end{align*}
By the local Lipschitz bounds of Corollary \ref{c:local Lipschitz bound} and the Maximum Principle, for sufficiently large $C(n, \gamma)$ in the definition of $h$, $u_{x, 0^+}(y) \le h(y)$ for all $y \in \{u_{x,0^+}>0 \} \cap B_1(0)$. 

We now do some quick geometric analysis on the auxiliary function $h$.  By standard harmonic theory, since $D_c$ is a cone in $B_1(0)$ the Almgren frequency function,
\begin{align*}
    N(0, r, h) = r\frac{\int_{B_r(0)}|\nabla h|^2 dx}{\int_{\partial B_r(0)}h^2 d\sigma}
\end{align*}
is a monotone non-decreasing function of $0< r < 1$. Furthermore the blow-up sequence $T_{0, r}h \rightarrow h_{\infty}$ strongly in $W^{1, 2}(B_1(0))$ as $r \rightarrow 0^+$ for a function $h_{\infty}$ which is homogeneous of degree $N(0, 0^+, h)= \lim_{r \rightarrow 0}N(0, r, h)$ and $\Delta h_{\infty} = 0$ in $\{ h_{\infty}>0\}$, the cone over $\partial B_1(0) \cap B_c(T_x\Gamma)$. See \cite{Mccurdy19} which proves these for star-shaped domains such as $D_c$.

We note that $u_{x, 0^+}(y) \le h(y)$ for all $y \in \{u_{x,0^+}>0 \} \cap B_1(0)$ it must be the case that $N(0, 0^+, h) \le \gamma+1$. We will show that for sufficiently small $0<c= c(n, \gamma)$ in the definition $D_c$, $N(0, 0^+, h) \ge \gamma +1.$

We begin with the symmetries of $\{h_{\infty}>0\} \subset \mathbb{R}^n$. Choose a unit normal $\vec v$ to the hyperplane $L$.  Since $\{h_{\infty}>0\}$ is invariant with respect to the rotations which preserve the hyperplane $L$ and is homogeneous, the function $h_{\infty}$ may be written as a function of two variables,
\begin{align*}
h_{\infty}(\rho, \theta) = R(\rho)T(\theta)  
\end{align*}
where $\rho = |\vec x|$ and $\theta \in [-\frac{\beta}{2}, \frac{\beta}{2}]$, is the angle that the vector $\vec x$ raises above $L$, with positive $\theta$ chosen in the $\vec v$ direction. We may write out the Laplace equation in terms of these variables,
\begin{align*}
   \Delta h & =  h_{\rho \rho} + \frac{n-1}{\rho}h_{\rho} + \frac{1}{\rho^2} \Delta_{\mathbb{S}^{n-1}}h \\
   & = h_{\rho \rho} + \frac{n-1}{\rho}h_{\rho} + \frac{1}{\rho^2}h_{\theta \theta} - \frac{(n-2)\tan(\theta)}{\rho^2}h_{\theta} = 0.
\end{align*}

Since $h_{\infty}(\rho, \theta) = R(\rho)T(\theta)$, we obtain
\begin{align*}
    0 & = R'' + (n-1)\frac{1}{\rho}R' - \lambda\frac{1}{\rho^2}R\\
    0 & = T'' - (n-2)\tan(\theta)T' - \lambda T.
\end{align*}
for some $\lambda \in \mathbb{R}$. We note that $R(\rho) = \rho^{c(\lambda)}$ for $0< c(\lambda) = \frac{-(n-2)+ \sqrt{(n-2)+4\lambda}}{2} < \infty$.  This forces $0< \lambda< \infty$ and shows that $c(\lambda) \rightarrow \infty$ as $\lambda \rightarrow \infty.$ We now consider the function $T$.  By Steiner symmetrization with respect to the hyperplane $L$ we see that we may assume that for $0<\beta$ small enough, $T'(\theta) \le 0$ for $0 \le \theta \le \frac{\beta}{2}.$ Therefore, if we consider a scalar multiple of $h_{\infty}$ (which we relabel as $h_{\infty}$) such that $\max\{h_{\infty}(y): y \in \partial B_{1}(0)\}=1$ it must be that $\max\{T(\theta): \theta \in [-\frac{\beta}{2}, \frac{\beta}{2}]\} = 1$.  Hence, we are able to estimate,
 \begin{align*}
     |T''(\theta)| \le \lambda
 \end{align*}
for all $0 \le \theta \le \frac{\beta}{2}$.  
 
By a second order Taylor expansion without remainder, for any $\theta \in [0, \beta/2]$
$$
T(\theta)= 1 + T'(0)(\theta) + \frac{1}{2}T''(\theta')\theta^2
$$
for some $\theta' \in [0, \theta]$. Since $T(\frac{\beta}{2}) = 0$ and $T'(0) = 0$ there must be a point $0 \le \theta' \le \frac{\beta}{2}$ such that $\frac{8}{\beta^2} = |T''(\theta')| \le \lambda$.  Therefore, as $c \rightarrow 0$ in the definition of $D_c$, $\beta \rightarrow 0$ and $\lambda \rightarrow \infty$, which implies that the degree of homogeneity of $h_{\infty}$ $c(\lambda) \rightarrow \infty$, as well.  However, this contradicts assumption that $c(\lambda) \le 1 + \gamma.$ Thus, we see that 
\begin{align*}
    c(n, \gamma) &\le \sin\left(\frac{\sqrt{2}}{\sqrt{(1+\gamma)(n-2)+ (1+\gamma)^2}}\right) \text{ implies   } 
    \lambda \ge \frac{8}{(1+\gamma)(n-2)+ (1+\gamma)^2}
\end{align*}
and $c(\lambda) \ge 1 + \gamma$. This proves the lemma.
\end{proof}

\subsection{Corollaries of the density lower bound}

\begin{lemma}\label{l:S closed}
Let $\Gamma$ be a $k$-dimensional $C^{1, \alpha}$-submanifold in $\mathbb{R}^n$, $0 \in \Gamma$, and let $Q(x) = \dist(x, \Gamma)^{\gamma}$. Let $u$ be a local minimizer of $J_{Q}(\cdot, B_2(0))$, then the set of  non-degenerate singular points $\partial \{u>0\} \cap \Gamma \setminus \Sigma$ is closed.
\end{lemma}

\begin{proof}
We note that $\Gamma \cap \partial \{u>0\}$ is closed.  The claim follows from the fact that the function 
\begin{align*}
x \mapsto W_{\gamma +1}(x, 0^+, u, \Gamma)
\end{align*}
is upper-semicontinuous in $\Gamma \cap \partial \{u>0 \} \cap B_1(0).$  Thus, if $\{x_i\}_i \subset \Gamma \cap \partial \{u>0\}$ then there is an $x \in \Gamma \cap \partial \{u>0\}$ and a subsequence $x_j$ such that $x_j \rightarrow x$ and $W_{\gamma +1}(x_j, 0^+, u, \Gamma) \le W^{\gamma +1}(x, 0^+, u)$.  By Lemma \ref{l:energy lower-bound}, $W_{\gamma+1}(x_j,0^+, u) \ge c_0>0$. Thus, $x \not \in \Sigma.$
\end{proof}

This next lemma allows us to find a scale at which the scalings $T_{x, r}u$ and $u_{x, r}$ are comparable.

\begin{lemma}\label{l:initial scale}
Let $\Gamma$ be a $k$-dimensional $(1, M)$-$C^{1, \alpha}$ submanifold with $0 \in \Gamma$. Let $\overline{B_2(0)} \subset \Omega$.  Let $u$ be an $\epsilon_0$-local minimizer of $J_{Q}(\cdot, \Omega)$ for $Q(x) = \dist(x, \Gamma)^{\gamma}$.  Suppose that $x \in \partial \{ u>0\} \cap \Gamma$.  There exists an $0<\eta_0(n, \gamma, \alpha)$ such that if $M \le \eta_0$, then there exists a $0<C_2(n, \gamma)$ such that if $x \in \mathcal{S}$, then for every scale $0< r< r_0$ and every $\overline{B_{2r}(x)} \subset \Omega,$
\begin{align*}
\dashint_{\partial B_1(0)} u^2_{x, r}(y)d\sigma(y) \ge C_2.
\end{align*}
In particular, then, there exist a pair of constants $0< c(n, \gamma) < C(n, \gamma) < \infty$ such that 
\begin{align}\label{e:non-degeneracy of rescaling}
    c(n, \gamma)u_{x, r}(u) \le T_{x, r}u(x) \le C(n, \gamma)u_{x, r}.
\end{align}
\end{lemma}

\begin{proof}
By Lemma \ref{l:energy lower-bound} and Corollary \ref{c:almost monotonicity}, if $[\Gamma]_{\alpha} \le M \le c_0(n, \gamma)\frac{1}{2C(n, \gamma, \alpha)} $, then we have that $W_{\gamma +1}(x, 0^+, u, \Gamma) \ge c_0$ and
\begin{align*}
    W_{\gamma + 1}(x, R, u) \ge \frac{1}{2}c_0,
\end{align*}
for all $0 \le R \le 1$. Therefore, 
\begin{align*}
    \frac{1}{2}c_0 & \le W_{\gamma +1}(x, R, u, \Gamma) \\
    & \le \frac{1}{R^{n+2\gamma}}\int_{B_R(x)}|\nabla u|^2 + Q^2(y)\chi_{\{u>0\}}dy\\
    & \le \frac{2C(n)\max_{y \in \{u>0\} \cap B_{R}(x)}\{\dist(y, \Gamma) \}^{2\gamma}\mathcal{H}^{n}(\{u>0\} \cap B_R(x))}{R^{n+2\gamma}}.
\end{align*}
Note, then, that if $\max_{y \in \{u>0\} \cap B_{R}(x)}\{\dist(y, \Gamma) \} \le CR$ for $C$ sufficiently small with respect to $c_0(n, \gamma)$, then we have a contradiction.  Thus, for all $0<R<1$, there is a point $y \in \{u>0 \} \cap B_R(x)$ such that $\dist(y, \Gamma) \ge C(n, \gamma) R$.  For example, $C(n, \gamma)$ may be taken to be a scalar multiple of $c(n, \gamma)$ in Lemma \ref{l:energy lower-bound}.

Now, if $B_{\frac{1}{2}C(n,\gamma)R}(y) \cap \partial \{u>0\} = \emptyset$, then Lemma \ref{l:interior balls} implies that $u(y) \ge C(n, \gamma)R^{1 + \gamma}.$  By the Maximum Principle, we can find a point $y_R \in \partial B_R(x)$ such that  $u(y_0) \ge C(n, \gamma)R^{1 + \gamma}$.  By the Lipschitz bound on $u$, there exists a ball of radius $C(n, \gamma)R$ in which  $u \ge \frac{1}{2}C(n, \gamma)R^{1+\gamma}.$ Integrating,
\begin{align*}
    \dashint_{\partial B_R(x)}u^2(y)d\sigma(y) \ge C(n, \gamma)R^{2(1 + \gamma)}.
\end{align*}

On the other hand, if $y_0 \in B_{\frac{1}{2}CR}(y) \cap \partial \{u>0\}$, then by Lemma \ref{l:qmin ineq},
\begin{align*}
    \dashint_{\partial B_{\frac{1}{4}C(n \gamma)R}(y_0)}u(y)d\sigma(y) \ge C_{\min}(\frac{1}{4}C(n, \gamma)R)^{1 + \gamma}.
\end{align*}
Thus, there must exist a point $z \in \partial B_{\frac{1}{4}C(n, \gamma)R}(y_0) \cap \{u>0\}$ such that $$u(z) \ge C_{\min}(\frac{1}{4}C(n, \gamma)R)^{1 + \gamma}.$$  Since by Lemma \ref{c:local Lipschitz bound}, $\text{Lip}(u|_{B_{\frac{1}{4}C(n, \gamma)R}}(y_0)) \le C(n)R^{\gamma}$, it must be the case that there exists a ball $B_{C(n, \gamma)R}(z)$ such that $u \ge C(n, \gamma)R^{1 + \gamma}$ on $B_{C(n, \gamma)R}(z).$  Thus, the desired lower bound holds,
\begin{align*}
    \dashint_{\partial B_R(x)}u^2(y)d\sigma(y) \ge c(n, \gamma)R^{2(1 + \gamma)}.
\end{align*}
The upper bound follows from Corollary \ref{c:local Lipschitz bound}.
\end{proof}

\begin{corollary}\label{c:Tu and u_r equivalence}
There is an $0< r_2(n, \gamma, \alpha, [\Gamma]_{\alpha}, \epsilon_0, \Lambda, A)$ such that the statement of Theorem \ref{t:compactness} holds with rescalings $T_{x, r}u$ replacing rescalings $u_{x, r}$ for all $0< r< r_2$. 
\end{corollary}

\begin{proof}
For the $\eta_0(n,\gamma, \alpha)>0$ the constant given by Lemma \ref{l:initial scale}, we let $0<r_1$ be a constant such that $[\Gamma]_{\alpha}r_1^{\gamma} \le \eta_0$.  Under this condition, note that if $x_0 \in \Gamma$, $[\Gamma^{x_0, r}]_{\alpha} \le \eta_{0}$ for all $0< r< r_1$.  We define,
\begin{align*}
    r_2 = \min\{r_1, r_0\}
\end{align*}
where $r_0$ is the standard scale for $\epsilon_0$-local minimizers in the class $K_{u_0, \Omega}$.
\end{proof}

\begin{remark}\label{r:non-degeneracy}(Non-degeneracy)
Let $\Gamma$ be a $k$-dimensional $C^{1, \alpha}$-submanifold in $\mathbb{R}^n$, $0 \in \Gamma$, and let $Q(x) = \dist(x, \Gamma)^{\gamma}$. If $\{u_i\}_i$ are a sequence of local minimizers of $J_{Q}(\cdot, B_2(0))$ with standard scale identically $r_0^i=1$, $x_i \in \mathcal{S}(u_i)$, $0<r_i \le r_2$ is a sequence of scales, then for any subsequential limit function $u_0$ obtained through Lemma \ref{t:compactness}, $u_0$ is non-degenerate in the sense that \eqref{e:non-degeneracy of rescaling} holds.
\end{remark}

\section{Proof of Theorem \ref{t:main theorem}}\label{s:proof of main theorem}

Let $u$ be an $\epsilon_0$-local minimizer in the class $K_{u_0, B_2(0)}$ with $\norm{\nabla u}^2_{L^2(B_2(0))} \le \Lambda$ and $\sup_{\partial B_{2}(0)} u_0 \le A$. Let $0< \epsilon$,  $0<\rho \le r$, and the integer $0 \le j \le k$ in $\mathcal{S}_{\epsilon, \rho}^j$ be given.  The strategy of the proof is to consider a scale $0<r_3$ at which we can employ the techniques of \cite{EdelenEngelstein19}. If $r_3 \le r$, we can easily cover $\mathcal{S}_{\epsilon, \rho}^j \cap B_1(0)$ by $C(n, k)r_3^{-k}$ balls of radius $r$.  If, on the other hand, $r \le r_3$, then, in each ball, we follow the proof of \cite{EdelenEngelstein19} Theorem 1.11, using Theorem \ref{t:compactness} and Lemma \ref{l:stong convergence} applied to the rescalings $T_{x, r}u$ in place of \cite{EdelenEngelstein19} Theorem 1.3, Lemma \ref{l:EE quant rigid} in place of \cite{EdelenEngelstein19} Lemma 3.1, and Lemma \ref{l:EE key dichot} in place of \cite{EdelenEngelstein19} Lemma 3.3.  The rest of the argument of \cite{EdelenEngelstein19} follows verbatim, save that where Edelen and Engelstein need to make the Holder seminorm $[Q]_{\alpha}$ small, we need to make the Holder seminorm $[\Gamma]_{\alpha}$ small. We define $0<r_3$ in Definition \ref{d: r3 def}.

For the benefit of the reader, we adumbrate the proof, below.  Following \cite{NaberValtorta17}, the proof falls broadly into two parts.  The first part consists of two Reifenberg-type results, the so-called Discrete Reifenberg and Rectifiable Reifenberg Theorems, which provide the necessary packing estimates.  These theorems are purely geometric measure theoretic results, which we restate, below.  We refer the interested reader to \cite{NaberValtorta17} for their proofs.

\begin{theorem}\label{t:discrete reifenberg}(Discrete Reifenberg \cite{NaberValtorta17})
Let $\{B_{r_q}(q)\}_q$ be a collection of disjoint balls with $q \in B_1(0)$ with radii $0< r_q \le 1$ and let $\mu$ be the packing measure $\mu = \sum_{q}r_q^{k}\delta_q$, for $\delta_q$ the Dirac measure at $q.$  There exist constants $\delta_{RD}(n) >0$ such that if 
\begin{align*}
    \int_{0}^{r} \int_{B_{2r}(x)} \beta^k_{\mu, 2}(z, s)^2d\mu \frac{ds}{s} \le \delta_{DR}r^k
\end{align*}
for all $0<r \le 1$ and all $x \in B_1(0)$, then $\mu(B_1(0)) = \sum_{q}r_q^k \le C_{DR}(n)$.
\end{theorem}

\begin{theorem}(Rectifiable Reifenberg \cite{NaberValtorta17})
Let $S \subset \mathbb{R}^n$.  There exists a constant $\delta_{RR}(n)>0$ such that if 
\begin{align*}
    \int_0^{2r}\int_{B_r(x)}\beta^k_{\mathcal{H}^k|_{S}, 2}(z, s)d\mathcal{H}^k|_{S}(z)\frac{ds}{s} \le \delta_{RR}r^k
\end{align*}
for all $x \in B_1(0)$ and $0<r\le 1$, then $S$ is countably $k$-rectifiable and $\mathcal{H}^k(S \cap B_r(x)) \le C_{RR}(n)r^k$ for all $x \in B_1(0)$ and $0<r\le 1$.
\end{theorem}

The second part of the argument consists of constructing a finite number of collections of balls to which we may apply these theorems.  Broadly speaking, the construction of these covers rely upon three observations and a complicated stopping-time argument.  Below, we state the three necessary observations, with careful details of any changes that need to be made from the proof in \cite{EdelenEngelstein19}.  However, we shall only sketch the stopping-time argument, as the details follow \cite{EdelenEngelstein19} verbatim.  

\subsection{Necessary Ingredients}
The observations are a kind of quantitative stability.  These type of results follow a predictable ``limit-compactness" argument by contradiction. 

Briefly, Lemma \ref{l:EE quant rigid} is a quantitative version of Lemma \ref{l:homogeneity}, wherein we showed that if $W_{\gamma +1}(x, r, u, \Gamma)$ is constant for $0< r \le 1$, then $u$ is $(\gamma + 1)$-homogeneous, and therefore is $(0, 0)$-symmetric in $B_1(x)$.  In \ref{l:EE quant rigid}, below, we show that if $W_{\gamma + 1}(x, r, u)$ almost constant, then $u$ is almost symmetric in $B_1(x).$

\begin{lemma}\label{l:EE quant rigid}(Quantitative Rigidity)
Fix $0 \le k \le n-1$.  Let $\delta>0$ and let $u$ be an $\epsilon_0$-local minimizer of $J_{Q}(\cdot , B_2(0))$ with $\Gamma$ a $k$-dimensional $C^{1, \alpha}$ submanifold.  Assume that $0< r_0 = 1$ and that $0 \in \Gamma \cap \partial \{u>0 \}$. Then there is a $\gamma_1 = \gamma(n, \delta)>0$ such that if $[\Gamma]_{\alpha} \le \gamma_1,$ and 
\begin{align*}
W_{\gamma + 1}(0, 1, u, \Gamma) - W_{\gamma + 1}(0, \gamma_1, u, \Gamma) \le \gamma_1,
\end{align*}
then $u$ is $(0, \delta)$-symmetric in $B_1(0)$.  
\end{lemma}

\begin{proof}
We argue by contradiction.  Suppose that there were a sequence of $\gamma_i \rightarrow 0$ and $u_i$ solutions to $J_{Q}(\cdot , B_2(0))$ with $\Gamma_i$ a $k$-dimensional $C^{1, \alpha}$ submanifold such that $0 \in \Gamma_i \cap \partial \{u_i>0 \}$, $[\Gamma_i]_{\alpha} \le \gamma_i$, and 
\begin{align*}
W_{\gamma + 1}(0, 1, u_i, \Gamma_i) - W_{\gamma + 1}(0, \gamma_i, u_i, \Gamma_i) \le \gamma_i,
\end{align*}
but $u_i$ is not $(0, \delta)$-symmetric in $B_1(0)$ for any $i \in \mathbb{N}.$
By Theorem \ref{t:compactness}, we may pass to a subsequence under which 
\begin{enumerate}
\item $\Gamma_j \rightarrow \Gamma$ a $k$-dimensional linear subspace.
\item $u_j \rightarrow u$ strongly in $W^{1, 2}(B_2(0))$ and in $C^{0, 1}(B_2(0)).$
\item By the aforementioned convergence, $W_{\gamma +1}(0, r, u_i, \Gamma_i) \rightarrow W_{\gamma + 1}(0, r, u, \Gamma)$.
\end{enumerate}

By Corollary \ref{c:limit solution}, the function $u$ satisfies the monotonicity formula, and by Corollary \ref{c:Weiss limit} we have that  
\begin{align*}
W_{\gamma + 1}(0, 1, u, \Gamma) - W_{\gamma + 1}(0, 0^+, u, \Gamma) & = \int_{0}^1 \frac{1}{s^{n + 1 + 2\gamma}}\int_{\partial B_s(0)}(\nabla u \cdot x - (\gamma +1)u)^2d\sigma(x) ds \\
& = 0.
\end{align*}
Hence, arguing as in Lemma \ref{l:homogeneity}, we have that $u$ is $(\gamma+ 1)$-homogeneous and is $(0, 0)$-symmetric by the non-degeneracy of Remark \ref{r:non-degeneracy}. Since $u_i \rightarrow u$ in $L^2(B_2(0)),$ this contradicts the assumption that each $u_i$ was not $(0, \delta)$-symmetric.
\end{proof}

The next observation is quite similar in nature. It is a quantitative version of the following fact: if $u$ is homogeneous with respect to the origin and translation invariant (i.e., symmetric) along a linear subspace $L$, and $u$ is homogeneous with respect to a point $x \not \in L$, then $u$ is symmetric with respect to $\text{span}\{x, L\}.$  In brief, the lemma below says that if $0 \in \mathcal{S}^k_{\epsilon}$ and $u$ is far from being $(k+1)$-symmetric in $B_1(0)$, but $u$ is also almost symmetric in $B_1(0),$ then $\mathcal{S}^{k}_{\epsilon}$ must be contained in the tubular neighborhood of a $k$-plane.  

\begin{lemma}\label{l:EE Beta sums}
Let $\epsilon>0$. Let $u$ be an $\epsilon_0$-local minimizer of $J_{Q}(\cdot, B_{10r}(x))$ for $0< 10r \le r_0$, $Q(x)=\dist(x, \Gamma)^\gamma$, and $\Gamma$ a $k$-dimensional $(1, \delta)$-$C^{1, \alpha}$ submanifold.  There is a $\delta(n, k, \gamma, \epsilon)>0$ such that if
\begin{align*}
\begin{cases}
u \text{  is $(0, \delta)$-symmetric in } B_{8r}(x)\\
u \text{  is  NOT $(k+1, \epsilon)$-symmetric in } B_{8r}(x),
\end{cases}
\end{align*}
then for any finite Borel measure $\mu$,
\begin{align*}
\beta^k_{\mu, 2}(x, r)^2 \le \frac{C(n, \alpha, \gamma, \epsilon)}{r^k}\int_{B_r(x)}W_{\gamma +1}(y, 8r, u, \Gamma) - W_{\gamma + 1}(y, r, u, \Gamma) + [\Gamma]_{\alpha}(8r)^\alpha d\mu(y).
\end{align*}
\end{lemma}

\begin{proof}
Fix the ball $B_r(x)$ and let $\mu$ be a finite Borel measure. We shall use the notation $A_{r, R}(x)$ to denote the annulus $B_{R}(x) \setminus B_r(x).$

This lemma follows from the analysis of the following non-negative quadratic form,
\begin{align*}
    B(v, w):= \dashint_{B_r(x)}(v \cdot (y- X))(w \cdot (y - X))d\mu(y)
\end{align*}
where $X = \dashint_{B_r(x)}yd\mu(y)$ is the center of mass. We note that $B$ has an orthogonal eigenbasis $v_1, ... , v_n$ with associated eigenvalues $\lambda_1 \ge \lambda_2 \ge ... \ge \lambda_n \ge 0$. Note that
\begin{align*}
    \beta^k_{\mu, 2}(x, r)^2 = \frac{\mu(B_r(x))}{r^k}(\lambda_{k+1} + ... + \lambda_n).
\end{align*}

The first estimate we need is that for sufficiently small $\delta>0$, and \textit{any} orthonormal basis, $\{v_i\}_i^{n}$, there exists a constant $0< C(n, \gamma, \alpha, \epsilon)< \infty$,
\begin{align*}
    \frac{1}{C(n, \gamma, \alpha, \epsilon)}\le r^{-n-2}\int_{A_{3r, 4r}(x)}\sum_{i=1}^{k+1}(v_i \cdot \nabla u(z))^2dz
\end{align*}
This follows from a straight-forward ``limit-compactness" argument, since for a $(0, 0)$-symmetric function if $\sum_{i=1}^{k+1}(v_i \cdot \nabla u(z))^2 = 0$ then $u$ is $(k+1, 0)$-symmetric if it is non-trivial. Note that non-triviality follows from non-degeneracy. The details follow \cite{EdelenEngelstein19} Theorem 5.1 verbatim, substituting Theorem \ref{t:compactness} and Theorem \ref{l:stong convergence} in place of \cite{EdelenEngelstein19} Theorem 1.3.

The second estimate requires greater detail.  We begin by observing that for any integers $0 \le i \le n$ and any $z$,
\begin{align*}
& \lambda_i(v_i \cdot \frac{1}{r^{\gamma}}\nabla u(z)) = B(v_i, \frac{1}{r^{\gamma}}\nabla u(z))\\
    &\qquad \qquad  = \dashint_{B_r(x)}(v_i \cdot (y - X)(\frac{1}{r^{\gamma}}\nabla u(z) \cdot (y + z - z - X)))d\mu(y) \\
    & \qquad \qquad \qquad \qquad - \frac{(\gamma +1)}{r^{\gamma}}u(z)\dashint_{B_r(x)}v_i \cdot (y - X)d\mu(y)\\
    & \qquad \qquad = \dashint_{B_r(x)}(v_i \cdot (y - X))(\frac{1}{r^{\gamma}}\nabla u(z) \cdot (y -z) -\frac{(\gamma +1)}{r^{\gamma}}u(z))d\mu(y)\\
    & \qquad \qquad \le  \lambda_i^{\frac{1}{2}}\left(\dashint_{B_r(x)}|\frac{1}{r^{\gamma}}\nabla u(z) \cdot (y -z) -\frac{(\gamma+1)}{r^{\gamma}}u(z)|^2 d\mu(y) \right)^{\frac{1}{2}}.
\end{align*}
Thus, $\lambda_i(v_i \cdot \frac{1}{r^{\gamma}}\nabla u(z))^2 \le \left(\frac{1}{r^{2\gamma}}\dashint_{B_r(x)}|\nabla u(z) \cdot (y -z) -(\gamma +1)u(z)|^2 d\mu(y) \right)^{\frac{1}{2}}.$  We continue calculating.
\begin{align*}
    &\lambda_i r^{-n-2}\int_{A_{3r, 4r}(x)}(v_i \cdot \nabla u(z))^2 dz\\
    & \qquad = r^{-n-2} \int_{A_{3r, 4r}(x)} \left(\frac{1}{r^{2\gamma}}\dashint_{B_r(x)}|\nabla u(z) \cdot (y -z) -(\gamma +1)u(z)|^2 d\mu(y) \right) dz\\
    &\qquad \le  5^{n+2+2\gamma}\dashint_{B_r(x)} \int_{A_{3r, 4r}(x)} \frac{1}{|z-y|^{n+2+2\gamma}}|\nabla u(z) \cdot (y -z) -(\gamma +1)u(z)|^2 dz d\mu(y)\\
    & \qquad = 5^{n+2+ 2\gamma}\dashint_{B_r(x)} \int_{A_{1r, 8r}(y)} \frac{1}{|z-y|^{n+2+2\gamma}}|\nabla u(z) \cdot (y -z) -(\gamma +1)u(z)|^2 dz d\mu(y)\\
    & \qquad = 5^{n+2}\dashint_{B_r(x)} \int_{1r}^{8r}\left( \frac{d}{d\rho}W_{\gamma + 1}(y, \rho, u, \Gamma) + 16\gamma [\Gamma]_{\alpha}r^{\alpha -1}\right) d\rho d\mu(y)\\
    & \qquad \le C(n, \gamma, \alpha) \dashint_{B_r(x)} W_{\gamma + 1}(y, 8r, u, \Gamma) - W_{\gamma +1}(y, r, u, \Gamma)+ [\Gamma]_{\alpha}r^{\alpha} d\mu(y).
\end{align*}
With these two estimates, we have,
\begin{align*}
    & \beta^k_{\mu, 2}(x, r)^2 \le \frac{\mu(B_r(x))}{r^k}n\lambda_{k+1}\\
    &  \le \frac{\mu(B_r(x))}{r^k}n C(n, \epsilon, \gamma, \alpha) \left( \sum_{i =1}^{k+1}\lambda_{i}\frac{1}{r^{n+2}}\int_{A_{3r, 4r}(x)}(v_i \cdot \nabla u(z))^2dz \right)\\
    & \le \frac{\mu(B_r(x))}{r^k} C(n, \gamma, \alpha, \epsilon, k) \left(\dashint_{B_r(x)} W_{\gamma + 1}(y, 8r, u, \Gamma) - W_{\gamma +1}(y, r, u, \Gamma)+ [\Gamma]_{\alpha}r^{\alpha} d\mu(y)  \right)\\
    &\le \frac{1}{r^k}C(n, \gamma, \alpha, \epsilon, k) \int_{B_r(x)} W_{\gamma + 1}(y, 8r, u, \Gamma) - W_{\gamma +1}(y, r, u, \Gamma)+ [\Gamma]_{\alpha}r^{\alpha} d\mu(y).
\end{align*}
This proves the lemma.
\end{proof}

The last crucial observation is a dichotomy which is a quantitative version of the geometry in the homogeneous case.  Simply stated, if $u$ is $k$-symmetric, but not $(k+1)$-symmetric, then $\partial \{u >0\}$ can be contained in either a $k$- or $(k-1)$-dimensional subspace.  

\begin{lemma}(Key Dichotomy)
Let $0 \le k \le n-1$ be an integer, and let $\Gamma$ by a $k$-dimensional $(1, M)$-$C^{1,\alpha}$ submanifold.  Let $u$ be an $\epsilon_0$-local minimizer of $J_{Q}(\cdot, B_4(0))$ such that $r_0 = 2$, $0 \in \partial \{u>0 \} \cap \Gamma$.  Let $0 <\rho<1$, $0< \gamma'<1$, and $0<\eta'<1$ be small, fixed constants. Let $\sup_{z \in B_1(0)}W_{\gamma +1}(z, 2, u, \Gamma) \le E \in (0, E_0]$ Then there exists a constant $\eta_1(n,  \gamma, \alpha, \epsilon, \rho, \gamma_1, \eta', E_0)<< \rho$ such that if $0<\eta \le \eta_1$ and $M \le \eta_1$ then either 
\begin{itemize}
    \item[i.] $W_{\gamma +1}(x, 2, u, \Gamma) \ge E - \eta'$ for all $x \in \mathcal{S}^k_{\epsilon} \cap B_1(0)$, or,
    \item[ii.] There exists a $(k-1)$-dimensional affine subspace $L^{k-1}$ such that $\{x \in B_1(0): W_{\gamma +1}(x, 2\eta, u, \Gamma) \ge E - \eta' \} \subset B_{\rho}(L^{k-1}).$
\end{itemize}
\end{lemma}

Note that we may take $E_0 = C(n)$ per the proof of Corollary \ref{c:almost monotonicity}.  The proof of this lemma is immediate from the following lemma.

\begin{lemma}\label{l:EE key dichot}
Fix $0\le k \le n-1.$  There are constants, $$\eta_1(n, \epsilon, \rho, \gamma_1, \eta', \alpha, \gamma)<< \rho \text{  and  } \beta(n, k, \eta', \rho, \gamma_1, \epsilon, \alpha)<1,$$ such that the following holds: Let $u$ be a minimizer to $J_{Q}(\cdot, B_4(0))$ with $\Gamma$ a $k$-dimensional $C^{1, \alpha}$ submanifold such that $0 \in \Gamma \cap \partial \{u>0 \}$ and $[\Gamma]_{\alpha}\le \Lambda.$ Assume that $$\sup_{x \in B_1(0)}W_{\gamma +1}(x, 2, u, \Gamma) \le E \in [0, c_n].$$

Suppose that $\eta \le \eta_1$ and $[\Gamma]_{\alpha} \le \eta$, and there are points $y_0, ... , y_{j} \in B_1(0) \cap \Gamma \cap \partial \{u>0\}$ satisfying,
\begin{align*}
y_i \not \in B_{\rho}(\langle y_0, ..., y_{i-1} \rangle)& \\
W_{\gamma +1}(y_i, 2\eta, u, \Gamma) \ge E - \eta, & \qquad  \forall i = 0, 1, ..., j.
\end{align*}
Then, writing $L= \langle y_0, ..., y_{j} \rangle$, for all $x \in B_{\beta}(L) \cap B_1(0),$ 
\begin{align*}
W_{\gamma +1}(x, \gamma_1 \rho, u, \Gamma) \ge E -\eta'
\end{align*}
and
\begin{align*}
\mathcal{S}^{j}_{\epsilon, \eta} \cap B_1(0) \subset B_{\beta}(L).
\end{align*}
\end{lemma}

\begin{proof}
There are two claims.  We argue both by contradiction.  Suppose that the first fails.  That is, suppose that there were a sequence of $u_i$ solutions to $J_{Q}(\cdot , B_2(0))$ with $\Gamma_i$ a $k$-dimensional $C^{1, \alpha}$ submanifold such that $0 \in \Gamma_i \cap \partial \{u_i>0 \}$, $[\Gamma_i]_{\alpha} \le \eta_i \rightarrow 0$, and collections $E_j, y_{ij}, \eta_j,$ and $\beta_j$ which satisfy the hypotheses.  Suppose that $\beta_j \rightarrow 0$ but that for each $j \in \mathbb{N}$ there is an $x_j \in B_{\beta_j}(L_j) \cap B_1(0)$ such that  
\begin{align*}
W_{\gamma + 1}(x_j, \gamma_i, u_i, \Gamma_i) \le E - \eta'.
\end{align*}
By Theorem \ref{t:compactness}, we may pass the a subsequence and obtain a function $u$ such that,
\begin{enumerate}
\item $u_j \rightarrow u$ strongly in $W^{1, 2}(B_4(0))$ and in $C^{0, 1}(B_4(0))$. 
\item $\Gamma_i \rightarrow \Gamma$ in the Hausdorff metric and $\Gamma$ is a $k$-dimensional linear subspace.
\item $E_j \rightarrow E$, $y_{ij}\rightarrow y_i$, $L_j \rightarrow L$, $x_j \rightarrow x \in \overline{B_1(0)} \cap L$.
\end{enumerate}
We note that since $\rho$ is fixed, the $y_i$ span $L$.  We note that by Corollary \ref{c:Weiss limit} and Corollary \ref{c:limit solution}, $u$ satisfies,
\begin{enumerate}
\item $\sup_{z \in B_1(0)}W_{\gamma + 1}(z, 2, u, \Gamma) \le E.$
\item $W_{\gamma + 1}(x, \gamma \rho, u, \Gamma) \le E - \eta'.$
\item $W_{\gamma + 1}(y_i, 0^+, u, \Gamma) \ge E$ for each $i = 0, 1, 2, ... , j.$
\end{enumerate} 

Since $\Gamma$ is flat, Lemma \ref{t:Weiss almost monotonicity} gives that the Weiss density is monotone increasing.  Thus, $u$ is $(\gamma + 1)$-homogeneous at $y_i$.  Therefore, $u$ is translation invariant along $L$ in $B_{1 + \delta} \subset \cup_{i}B_2(y_i)$ for some $\delta >0.$  In particular, then, $W_{\gamma + 1}(x, 0^+, u, \Gamma) = E,$ which is a contradiction.

We now argue the second claim.  That is, fix $\beta$, and assume that we have sequences $u_j, E_j, y_{ij}, L_j,$ which satisfy the hypotheses of the lemma, and $\eta_j \rightarrow 0$ such that for $j$ there exists some $x_j \in \mathcal{S}^{j}_{\epsilon, \eta_j}(u_j) \cap B_1(0) \setminus B_{\beta}(L_j).$  
By Theorem \ref{t:compactness}, we may pass to a subsequence such that $u_j, E_j, y_{ij}, L_j$ converge as above.  Again, we have that for some $\delta>0$, the limit function $u$ will be $k$-symmetric with respect to $L$ in $B_{1+\delta}(0).$  Since the limit point $x \in \overline{B_1(0)} \setminus B_{\beta}(L),$ any blow-up of $u$ at $x$ will be $(k+1)$-symmetric.  In particular, for some $r>0$, $u_j$ will be $(k+1, \epsilon)$-symmetric in $B_r(x).$  This is a contradiction.
\end{proof}

\begin{definition}\label{d: r3 def}(The scale $0<r_3$)
We define the base scale $0< r_3$ using the constants from the lemmata, above.  Let the $0<\delta$ in Lemma \ref{l:EE quant rigid} be the $0<\delta$ from Lemma \ref{l:EE Beta sums}. Let $0<\eta'$ from Lemma \ref{l:EE key dichot} be the $0<\gamma_1(n, \delta)$ from Lemma \ref{l:EE quant rigid}.  Let $0<r_2(n ,\gamma, \alpha, [\Gamma]_{\alpha}, \epsilon_0, \Lambda, A)$ be as in Corollary \ref{c:Tu and u_r equivalence}.  We define $0< r_3$ to satisfy the following two conditions. First, we let $0<r_3$ be a constant which satisfies
\begin{align*}
0<r_3 \le \frac{1}{20}r_2, 
\end{align*}
where $r_2$ is as in Corollary \ref{c:Tu and u_r equivalence}.  Second, let $0<r_3$ be sufficiently small so that in each ball $x \in \Gamma$, $[\Gamma^{x, r_3}]_{\alpha} \le \min\{\gamma_1, \delta, \eta_1\}$, where $\eta_1(n, \epsilon, \rho, \gamma_1, \eta', \alpha, \gamma)$ is as in Lemma \ref{l:EE key dichot}. This makes $0< r_3 = r_3(n, k, \gamma, \alpha, [\Gamma]_{\alpha}, \epsilon_0, \Lambda, A, \epsilon, j)$.
\end{definition}

With these observations in hand, the argument becomes procedural, and the construction of the desired coverings follows the proof in Edelen-Engelstein verbatim, with Theorem \ref{t:compactness} and Lemma \ref{l:stong convergence} in place of \cite{EdelenEngelstein19} Theorem 1.3, Lemma \ref{l:EE quant rigid} in place of \cite{EdelenEngelstein19} Lemma 3.1, and Lemma \ref{l:EE key dichot} in place of \cite{EdelenEngelstein19} Lemma 3.3. 

In particular, these observations and Theorem \ref{t:discrete reifenberg} allow us to prove the following lemma.

\begin{lemma}\label{l:packing estimate}
Let $0 \le k \le n-1$ and $0< \gamma, \alpha$ be fixed.  Let $\Gamma$ be a $k$-dimensional $(1, M)$-$C^{1, \alpha}$ submanifold such that $0 \in \Gamma$, and let $Q(x) = \dist(x, \Gamma)^{\gamma}$. Let $u$ be an $\epsilon_0$-local minimizer for $J_{Q}(\cdot, B_3(0))$ with $0< r_0 = 2$, and let $E = \sup_{y \in B_1(0)}W_{\gamma + 1}(u, 2, y, \Gamma).$

Let $0< R.$  There is an $\eta (n, \gamma, \alpha, \epsilon, M_0) > 0$ such that if $[\Gamma]_{\alpha} \le \eta$ and $\{B_{\eta r_p}(p)\}_p$ is a collection of disjoint balls satisfying,
\begin{itemize}
    \item [i.] $W_{\gamma + 1}(u, \eta r_p, p, \Gamma) \ge E - \eta.$
    \item [ii.] $p \in \mathcal{S}^{j}_{\epsilon, R}.$
    \item [iii.] $R \le r_p \le 1.$
\end{itemize}
then $\{B_{\eta r_p}(p)\}_p$ satisfies the packing estimate,
\begin{align*}
    \sum_p r_p^j \le C(n).
\end{align*}
\end{lemma}

The proof follows that of \cite{EdelenEngelstein19} Lemma 6.1 \textit{mutatis mutandis}. The rest of the proof of Theorem \ref{t:main theorem} relies upon constructing a finite sequence collections of balls whose union covers each of the quantitative strata $\mathcal{S}_{\epsilon, R}^j$.

\subsection{Cover Construction}

The idea is to use the dichotomy implied by Lemma \ref{l:EE key dichot} to set up a stopping-time argument.  For ease of notation, we rescale.  In any ball $B_{r_3}(x)$, we dilate everything to the ball $B_1(0)$, where $\Gamma^{x, r_3}$ satisfies the conditions of the previous lemmata.  The dichotomy implies that either  $\mathcal{S}^{j}_{\epsilon, r} \cap B_1(0) \subset B_{\beta}(L) \cap B_1(0)$ for some $j$-dimensional affine linear subspace and the Weiss density does \textit{not} have a large drop on this set, or, the drop in the Weiss density must be \textit{large} for all points outside of a neighborhood of a lower-dimensional affine subspace $B_\rho(L^{j-1})$.  In the former case, we call $B_1(0)$ a ``good" ball.  In the latter case, we call $B_1(0)$ a ``bad" ball.

In a ``good" ball, the procedure is to cover $B_{\beta}(L) \cap B_1(0)$ by balls of a smaller scale and then use the dichotomy to inductively refine a cover, stopping in ``bad" balls, and refining the cover in ``good" balls.  The full details of the construction are in Section 7.1 of \cite{EdelenEngelstein19}. By construction, the collection of ``good" balls and ``bad" balls so produced satisfy the hypotheses of the Lemma \ref{l:packing estimate}.  

In a ``bad" ball, the construction becomes one step more complicated.  In a ``bad" ball, we again cover the set of points with \text{small} drop in Weiss density in the neighborhood $B_{\rho}(L^{j-1}) \cap B_1(0)$ by balls of a smaller scale, and we cover the rest of the ball in $B_1(0) \setminus B_{\beta}(L^{j-1})$, as well. We sort the balls covering $B_{\rho}(L^{j-1}) \cap B_1(0)$ into ``good" and ``bad" balls and  inductively refine in ``bad" balls.  We stop in ``good" balls.  This creates three collections of balls.  First, the collection of "bad" balls stemming from the covers of the tubular neighborhoods of the  $(j-1)$-dimensional subspaces.  Second, the ``good" balls stemming from the covers of the tubular neighborhoods of the $(j-1)$-dimensional subspaces. Third, the balls which cover the complement, $B_{r_q}(q) \setminus B_{\beta r_q}(L^{j-1}_{q})$ in the collection of``bad" balls. The full details of the construction are in Section 7.2 of \cite{EdelenEngelstein19}.

For the collections of ``good" and ``bad" balls in the ``bad" ball construction, the lower-dimensional containment in neighborhoods of affine $(j -1)$-dimensional subspaces easily give $j$-dimensional packing estimates which are stronger than Lemma \ref{l:packing estimate} (see Theorem 7.6 \cite{EdelenEngelstein19} for details).  

While there is no way to obtain $j$-dimensional packing estimates on the third collection, we do know that the Weiss density has large \textit{definite} drop in each of these balls. Hence, though we restart the ``good" tree or ``bad" tree construction in each of these balls, since the Weiss density is of bounded variation there is a fixed constant $C(n, \gamma, \alpha, \epsilon)<\infty$ such that we can restart this procedure at most $C(n, \gamma, \alpha, \epsilon)$ times. Thus, there are a finite number of such collections of balls, each of which satisfies the desires packing estimates.

The rectifiability of $\mathcal{S}^j_{\epsilon}$ follows the proof of \cite{EdelenEngelstein19} Theorem 1.11 \textit{verbatim}.

\section{Containment and Density upper bound}\label{s:e reg}

In this section we focus upon containment results.  Lemma \ref{l:containment minimizers}, below, proves Lemma \ref{l:e reg}.  Note that it relies upon the non-degeneracy of solutions.

\begin{lemma}\label{l:containment minimizers}(Containment for Local Minimizers)
Let $\Gamma$ be a $k$-dimensional $(1, \frac{1}{4})$-$C^{1, \alpha}$ submanifold such that $0 \in \Sigma$. If $Q(x) = \dist(x, \Gamma)^{\gamma}$ and $u$ is an $\epsilon_0$-local minimizer of $J_Q(\cdot, B_2(0))$ in the class $K_{u_0, B_2(0)}$ for $\norm{\nabla u_0}^2_{L^2(B_2(0))} \le \Lambda$ and $\sup_{y \in \partial B_2(0)}u_0(y) \le A$, then, if $k = n-1$ there exists an $0< \epsilon(n, \gamma, \alpha, M, \Lambda, A, \epsilon_0)$ such that $\mathcal{S} \subset \mathcal{S}^{n-2}_{\epsilon}.$
\end{lemma}

\begin{proof}
We argue by contradiction. Fix $\alpha > 0$, $\gamma>0$, and a positive number $M$. Suppose that for every $i \in \mathbb{N}$, there exists a $n-1$-dimensional $C^{1, \alpha}$ submanifold $\Gamma_i$ such that $0 \in \Gamma_i$ and $[\Gamma_i]_{\alpha} \le \frac{1}{4}$, a function $u^{(i)}$ an $\epsilon_0$-local minimizer of $J_{Q_i}(\cdot, B_2(0))$ , where $Q_i = \dist(x, \Gamma_i)^{\gamma}$, and a point $x_i \in \overline{B_1(0)} \cap \Gamma_i$ and a radius $0<r_i \le \dist(x_i, \partial B_2(0)) \le 2$ such that $x_i \in \partial\{u_i >0\}$ and there exists a normalized non-trivial $(n-1)$-symmetric function $\phi_i$ such that $\norm{u_{x_i, r_i} - \phi_i}_{L^2(B_1(0))} \le 2^{-i}.$  Then Theorem \ref{t:compactness} and Lemma \ref{l:manifold convergence} imply that there is a $(n-1)$-dimensional $C^{1, \alpha}$ submanifold $\Gamma_0$ such that $[\Gamma_0]_{\alpha} \le \frac{1}{4}$ and a function $u \in W^{1, 2}(B_1(0))$ such that, by passing to a subsequence,
\begin{align*}
    \Gamma_i \cap B_2(0) & \rightarrow \Gamma_0 \cap B_2(0)  \qquad \text{ in the Hausdorff metric}\\
    u^{(i)}_{x_i, \min\{r_i,r_0\}} & \rightarrow  u \text{  in $W^{1, 2}(B_1(0))$ and in $C^{0, 1}(B_1(0))$}.
\end{align*}
Furthermore, $u$ is $(n-1, 0)$-symmetric in $B_1(0)$ because $u$ is non-degenerate in the sense of Remark \ref{r:non-degeneracy}. Therefore, $u$ must be a non-trivial piece-wise linear function.  In particular, it must be $1$-homogeneous.  However, since by Corollary \ref{c:local Lipschitz bound} for all $y \in B_1(0)$
\begin{align*}
    |\nabla u^{(i)}_{x_i, \min\{r_i,r_0\}}(y)| \le C(n)2^{\gamma}|y|^{1+\gamma}. 
\end{align*}
Therefore, for a radius $0< \rho$ depending only upon $n$ and $\gamma$ there is a ball $B_\rho(0)$ such that $|\nabla u^{(i)}_{x_i, \min\{r_i,r_0\}}(y)| \le \frac{1}{2}c_n$ for all $y \in B_\rho(0)$, where $c_n = |\nabla u|$ in $\{u>0\}$. Since this holds for all $i$, we reach a contradiction with strong convergence in $W^{1, 2}(B_1(0))$.
\end{proof}

The upper bound on the $Q$-density $W_{1+\gamma}(x, 0^+, u, \Gamma)$ is a result of some weak results on $\{u=0\}$ near $\partial\{u>0\}$.  Stronger results than those proved below are known for local minimizers.  See, for example \cite{DavidToro15} in which an Exterior Ball Condition is proved for almost-minimizers. However, nothing so powerful is needed here.

\begin{lemma}\label{exterior mass condition}(cf. \cite{AltCaffarelli81} Lemma 3.7)
Let $u$ be as in the hypotheses of Lemma \ref{l:sigma and S decomp}. And let $0< r_0$ be the standard scale. For all $0< r< r_0$ and all $x \in \partial \{u>0\}$, if $\dist(x, \Gamma) \ge 2r$ then there is a constant $0<C(n, \gamma)$ such that 
\begin{align*}
    C(n, \gamma) r (\min_{w \in B_r(x)}\{Q(w) \})^2 \omega_n r^{n} \le  \int_{B_r(x)}Q^2 \chi_{\{u>0 \}}dx.
\end{align*}
\end{lemma}

\begin{proof}
We note that for $h^{x, r}_u$ the harmonic extension of $u$ in $B_R(x)$ we have from comparison and the Poincare inequality
\begin{align*}
r^{-1}\int_{B_r(x)}(u -h^{x, r}_u)^2 dx \le \int_{B_r(x)}Q^2\chi_{\{u=0\}}dx.
\end{align*}
By Corollary \ref{c:local Lipschitz bound}, $|\nabla u|, |\nabla h^{x, r}_u| \le C(n)\max_{w \in B_r(x)}\{Q(w)\}$ and by Corollary \ref{c:u growth} and Corollary \ref{l:interior balls} $$|h^{x, r}_u(y)| \ge C(n, \gamma) r \min_{w \in B_r(x)}\{Q(w)\}$$ for a set of $y \in \partial B_r(x)$ with measure $(cr/2)^{n-1}$ for 
$$
c = \frac{C_{min}\min_{w \in B_{\frac{1}{2}r}(x)}\{ Q(w)\}}{4C(n)\max_{B_{\frac{3}{4}r}(x)}\{Q(w)\}}.
$$
Thus, $|h^{x, r}_u(x)| \ge \frac{c}{2^{n-1}}C(n, \gamma) r \min_{w \in B_r(x)}\{Q(w)\}$.  Therefore there is a ball of radius $\rho \ge 2^{-n-1}c(n, \gamma) (\frac{\min_{w \in B_{r}(x)}\{ Q(w)\}}{\max_{B_{r}(x)}\{Q(w)\}})^2 r$ in which $|u-h^{x, r}_u| \ge \frac{c}{2^{n}}C(n, \gamma) r \min_{w \in B_r(x)}\{Q(w)\}$.  The claim follows immediately recalling $(\frac{\min_{w \in B_{r}(x)}\{ Q(w)\}}{\max_{B_{r}(x)}\{Q(w)\}})^2 \le 2^{2\gamma}$.
\end{proof}

\begin{lemma}\label{Q-density upper bound}(Proof of the $Q$-density upper bound)
For $u$ a local minimizer satisfying the hypotheses of Lemma \ref{l:sigma and S decomp}.  Then there is a $C(n, \gamma)< \int_{B_1(0)}|x_n|^{2\gamma}dx$ such that $\Phi(x, 0^+) \le C(n, \gamma)$ for all $x \in \Gamma \cap \partial \{u>0\}$.
\end{lemma}

\begin{proof}
We argue by contradiction.  Suppose that there were a sequence of $u_i,x_i$ as in the hypotheses for which $\Phi(x_i, 0^+)(u_i) \rightarrow \int_{B_1(0)}|x_n|^{2\gamma}dx.$  Letting $\tilde u_i$ be a blow-up of $u_i$ we may use Theorem \ref{t:compactness} to pass to a subsequence which converges to a non-trivial $(1+\gamma)$-homogeneous function $u$ for which $\int_{B_1(0)}Q^2 \chi_{\{u>0\}}dx = \int_{B_1(0)}|x_n|^{2\gamma}dx$.  That is, $\{u>0\}$ has full measure.  There are two cases.  First, if $\partial \{u>0\} \cap \partial B_1(0) = \emptyset$ then $u$ is a function which is harmonic in $\mathbb{R}^n \setminus \{0\}$ which is continuous at $\{0\}$. Therefore, $\{0\}$ is removable and $u$ is a non-negative enitre harmonic function with minimum acheived at $u(0)=0$.  This is absurd.  The second case is that $\partial \{u>0\} \cap \partial B_1(0) \not = \emptyset.$  In this case, since $x_i \in \partial \{\tilde{u}_i >0\} \cap \Gamma$ and $\partial \{\tilde{u}_i >0\}$ converges locally to $\partial \{u>0\}$ away from $\Gamma$ the exterior mass condition in Lemma \ref{exterior mass condition} must be violated for sufficiently large $i$.  This contradiction proves the claim.
\end{proof}

\section{Finite perimeter}\label{s:finite perimeter}

In this section, we focus upon the set of cusps $\Sigma = \partial \{u>0\} \cap \Gamma \setminus \mathcal{S}^{n-2}_{\epsilon_0}$.  
The main result of this section is the following Lemma.

\begin{lemma}\label{l:finite perimeter}
Let $\Gamma$ be a $k$-dimensional $(1, \frac{1}{4})$-$C^{1, \alpha}$ submanifold such that $0 \in \Sigma$. If $Q(x) = \dist(x, \Gamma)^{\gamma}$ and $u$ is a local minimizer of $J_Q(\cdot, B_2(0)$, then $\partial \{u>0 \} \cap B_1(0)$ is a set of finite perimeter.  In particular, if $k= n-1$, then set of cusp points $$\Sigma = \Gamma \cap \partial \{u>0\} \setminus \mathcal{S}$$ has $\mathcal{H}^{n-1}$-measure zero.
\end{lemma}

Before proving Lemma \ref{l:finite perimeter}, we need the results of the following lemmata.

\begin{lemma}\label{l:finite perimeter 1}(Inwards-Minimality implies Finite Perimeter Estimate)
Let $n \ge 2$ and fix $0< \gamma$, $0 < \alpha \le 1$, and an integer $0 \le k \le n-1$.  Let $\Gamma$ be a $k$-dimensional $(1, \frac{1}{4})$ $C^{1, \alpha}$-submanifold.  Let $Q(x) = \dist(x, \Gamma)^{\gamma}$ and let $u$ be an $\epsilon_0$-local minimizer of $J_{Q}(\cdot, B_2(0))$ in the class $K_{u_0, B_2(0)}$ for $\norm{\nabla u_0}^2_{L^2(B_2(0))} \le \Lambda$ and $\sup_{y \in \partial B_2(0)} u_0(y) \le A$.  Let $0< r_0$ be the standard scale.  Then, for all $0< r \le \frac{1}{2}r_0$ there is a constant $0<C<\infty$ such that
\begin{align*}\int_{\{0 < u \le \epsilon \} \cap B_{r}(x_0)}|\nabla u|^2 + Q^2(x)\chi_{\{0<u \le \epsilon \}}dx \le C \epsilon,\end{align*}
for all $\epsilon \in (0, 1]$. Moreover, we may take
\begin{align*}C = C(n)(||\nabla u||_{L^{2}(B_{2r}(x_0))}r^{\frac{n-2}{2}} + \epsilon r^{n-2}).
\end{align*}
\end{lemma}

\begin{proof}
Let $\phi \in C^{\infty}(\mathbb{R}^n)$ such that
\begin{align*}\phi(x) = \begin{cases}0 & x \in B_{r}(x_0)\\
1 & x \in \mathbb{R}^n \setminus B_{2r}(x_0)\end{cases}
\end{align*}
Note that we may choose $\phi$ to be radial, $\phi = \phi(|x- x_0|)$, such that $\phi > 0$ in $B_{2r}(x_0) \setminus B_r(x_0)$, $\phi < 1$ in $B_{2r}(x_0)$, and such that $|\nabla \phi| \le \frac{5}{r}.$

For a fixed $\epsilon,$ let $u_{\epsilon} = (u - \epsilon)_+$ and $\tilde{u}_{\epsilon} = \phi u + (1- \phi)u_{\epsilon}.$ Note that $\tilde u_{\epsilon} > 0$ in $B_{2r}(x_0) \setminus B_r(x_0)$ if $u > 0$, and $\tilde u_{\epsilon} > 0$ in $B_{r}(x_0)$ if and only if $u > \epsilon.$
\begin{align*}
|\nabla \tilde{u}_{\epsilon}|^{2} = & \chi_{\{0 < u \le \epsilon\}}|\nabla (u\phi)|^2 + \chi_{\{u>\epsilon\}}|\nabla( u - \epsilon(1-\phi))|^2\\
\le & \chi_{\{0 < u \le \epsilon\}}(\phi^2 |\nabla u|^2 + \epsilon 2 |\nabla u| |\nabla \phi| + \epsilon^2|\nabla \phi|^2) \\
& + \chi_{\{u>\epsilon\}}( |\nabla u|^2 + \epsilon 2 |\nabla u||\nabla \phi| + \epsilon^2 |\nabla \phi |^2)
\end{align*}
Then, if we let 
\begin{align*}
C & = C(n)(||\nabla u||_{L^{2}(B_{2r}(x_0))}||\nabla \phi||_{L^{2}(B_r(x_0))} + \epsilon ||\nabla \phi||^2_{L^2(B_{2r}(x_0))})\\
& \le C(n)(||\nabla u||_{L^{2}(B_{2r}(x_0))}r^{\frac{n-2}{2}} + \epsilon r^{n-2}),    
\end{align*}
then because $2r \le r_0$, by the local minimality of $u$ we may compare
\begin{align*}
0 \ge & \left(\int_{B_{2r}(x_0)} |\nabla u|^2 + Q^2\chi_{\{u> 0\}}dx \right) - 
\left( \int_{B_{2r}(x_0)} |\nabla \tilde{u}_{\epsilon}|^2 + Q^2\chi_{\{\tilde{u}_{\epsilon}> 0\}}dx \right)\\
\ge & \left(\int_{B_{2r}(x_0)} |\nabla u|^2 - \int_{B_{2r}(x_0)} |\nabla \tilde{u}_{\epsilon}|^2 dx \right) + \int_{B_{r}(x_0)} Q^2\chi_{\{0< u \le \epsilon \}}dx \\
\ge & \int_{B_{2r}(x_0) \cap \{0 < u \le \epsilon\}} (1-\phi^2) |\nabla u|^2 dx + \int_{B_{r}(x_0)} Q^2\chi_{\{0< u \le \epsilon \}}dx - C \epsilon \\
\ge & \int_{B_{r}(x_0) \cap \{0 < u \le \epsilon\}} |\nabla u|^2 dx + \int_{B_{r}(x_0)} Q^2\chi_{\{0< u \le \epsilon \}}dx - C \epsilon.
\end{align*}
\end{proof}

\begin{lemma}\label{l:finite perimeter 2}(Locally finite perimeter bounds)
Suppose that $D \subset \mathbb{R}^n$ is an open, bounded set, $Q \in C^{0, \alpha}(D)$, $0< Q_{D, min} \le Q$, and that $\phi: D \rightarrow [0, +\infty]$ is a function in $W^{1, 2}(D)$ such that there exists an $\overline{\epsilon} >0$ and a $C > 0$ such that for all $0< \epsilon \le \overline{\epsilon}$
\begin{align*}
\int_{\{0< \phi \le \epsilon \} \cap D} |\nabla \phi|^2 + Q^2(x)\chi_{\{0< \phi \le \epsilon \}}dx \le C \epsilon.
\end{align*}
Then, $Per(\{\phi > 0 \}; D) \le C\frac{1}{\sqrt{Q_{D, min}}},$ where $Q_{D, min} = \min_{y \in D}\{Q(y)\}$.  Furthermore, if $D = B_r(x)$, then $C$ may be taken to be $C_n\norm{\nabla u}_{L^2(B_{2r}(x))}r^{\frac{n-2}{2}}$ 
\end{lemma}
\begin{proof}
\begin{align*}\int_{0}^{\epsilon} \mathcal{H}^{n-1}(\{\phi = t \} \cap D)dt & = \int_{\{0 < \phi \le \epsilon \} \cap D} |\nabla \phi| dx\\& \le |\{0 < \phi \le \epsilon \} \cap D|^{\frac{1}{2}} (\int_{\{0 < \phi \le \epsilon \} \cap D} |\nabla \phi |^2dx)^{\frac{1}{2}}.\end{align*}
We now bound each term by our assumption.
\begin{align*}
\int_{0}^{\epsilon} \mathcal{H}^{n-1}(\{\phi = t \} \cap D)dt & \le |\{0 < \phi \le \epsilon \} \cap D|^{\frac{1}{2}} (\int_{\{0 < \phi \le \epsilon \} \cap D} |\nabla \phi |^2dx)^{\frac{1}{2}}\\
& \le \frac{Q_{D, min}}{Q_{D, min}}|\{0 < \phi \le \epsilon \} \cap D|^{\frac{1}{2}} (\int_{\{0 < \phi \le \epsilon \} \cap D} |\nabla \phi |^2dx)^{\frac{1}{2}}\\
& \le (C \epsilon \frac{1}{Q_{D, min}})^{\frac{1}{2}}(C \epsilon)^{\frac{1}{2}}\\
& \le C \epsilon \frac{1}{\sqrt{Q_{D, min}}}.
\end{align*}
Now, let $0<\epsilon$ and let $\delta \in (0, \epsilon]$ such that
\begin{align*}
\mathcal{H}^{n-1}(\partial^*(\{\phi > \delta\} \cap D)) & \le \frac{1}{\epsilon} \int_{0}^{\epsilon} \mathcal{H}^{n-1}(\{\phi = t \} \cap D) dt\\
& \le C \frac{1}{\sqrt{Q_{D, min}}}.
\end{align*} 
Then, letting $\epsilon \rightarrow 0$, $\delta \rightarrow 0$ and we obtain, $\mathcal{H}^{n-1}(\partial^*(\{\phi > 0\} \cap D)) \le C \frac{1}{\sqrt{Q_{D, min}}}.$
\end{proof}

\subsection{Proof of Lemma \ref{l:finite perimeter}}
Following Lemma \ref{l:distance estimate}, we decompose $B_{1}(0)$ into dyadic annular neighborhoods,
\begin{align*}
    A_{j} = \left(B_{2^{-j}}(\Gamma) \setminus B_{2^{-j-1}}(\Gamma)\right) \cap B_{1}(0).
\end{align*}
Now, using the estimate in Lemma \ref{l:distance estimate}, we can cover each $A_{j}$ by $C(n)2^{jk}$ balls of radius $2^{-j-3}$.  We denote this collection by $\mathcal{B}_j = \{B_j\}$.  For a ball $B_j$, we shall use $2B_j$ to denote the concentric dilate of $B_j$ by inflation factor $2$. By Lemma \ref{l:finite perimeter 1}, Lemma \ref{l:finite perimeter 2}, and the Lipschitz estimate from Corollary \ref{c:local Lipschitz bound} we obtain that for sufficiently large $j$ depending upon the ``standard scale" $0<r_0$, 
\begin{align*}
    \mathcal{H}^{n-1}(\partial \{u>0\} \cap B_j) & \le \frac{1}{\sqrt{Q_{\min, B_j}}} C(j)\\
    & \le C(j) (2^{-j-3})^{-\frac{1}{2}\gamma},
\end{align*}
where $C(j)$ is as in Lemma \ref{l:finite perimeter 1}
\begin{align*}
    C(j)  \le C(n)(||\nabla u||_{L^{2}(2B_j)}(2^{-j-3})^{\frac{n-2}{2}}) \le C(n)(C(n)2^{\gamma} (2^{-j})^{\frac{n + 2\gamma}{2}}(2^{-j-3})^{\frac{n-2}{2}}).
\end{align*}
Note that since we took $\epsilon \rightarrow 0$ in Lemma \ref{l:finite perimeter 2}, the $\epsilon$ term in $C$ from Lemma \ref{l:finite perimeter 1} vanished.

Summing over each of the $C(n)2^{jk}$ in the collection $\mathcal{B}_j$, we obtain the estimate,
\begin{align*}
    \mathcal{H}^{n-1}(\partial \{u>0\} \cap A_{j}) \le C(n)2^{-j(n-k-1+\gamma) -3\frac{n-2}{2} + \frac{5}{2}\gamma}.
\end{align*}
Summing over $j$ and recalling that $k \le n-1$, we obtain that
\begin{align*}
    \mathcal{H}^{n-1}(\partial \{u>0 \} \cap B_{2^{j_0}}(0)) & = \mathcal{H}^{n-1}(\Gamma) + \mathcal{H}^{n-1}(\partial \{u>0 \} \cap \bigcup_{j= j_0}^{\infty}A_{j})\\
    & \le C(n, [\Gamma]_{\alpha})\omega_{n-1} + C(n) \sum_{j = 1}^{\infty}2^{-j(n-k-1+\gamma) -3\frac{n-2}{2} + \frac{5}{2}\gamma}\\
    & \le C(n, [\Gamma]_{\alpha})\omega_{n-1} + C(n, \gamma) \sum_{j = 1}^{\infty}2^{-j(\gamma)} \le C(n, [\Gamma]_{\alpha}, \gamma) <\infty.
\end{align*}

Since $\partial \{u>0\} \cap B_1(0) \setminus B_{r^{-j -3}}(\Gamma)$ has finite $\mathcal{H}^{n-1}$-measure for every integer $j$, and the previous estimate shows $\partial \{u>0\} \cap B_1(0) \cap B_{r^{-j_0 -3}}(\Gamma)$ has finite $\mathcal{H}^{n-1}$-measure for $j_0$ sufficiently large, $\{u>0\}$ is a set of locally finite perimeter, and hence $\text{sing}(\partial \{u>0\})$ has $\mathcal{H}^{n-1}$-measure zero.

\bibliography{references}
\bibliographystyle{amsalpha}

\end{document}